\documentclass{amsart}
\usepackage{amssymb}
\usepackage[all]{xy}

\newcommand{\C}{{\mathbb{C}}}

\newcommand{\N}{{\mathbb{N}}}
\newcommand{\Q}{{\mathbb{Q}}}
\newcommand{\R}{{\mathbb{R}}}

\newcommand{\Z}{{\mathbb{Z}}}

\newcommand{\Ch}{{\mathcal C}}
\newcommand{\Dh}{{\mathcal D}}
\newcommand{\Eh}{{\mathcal E}}

\newcommand{\Kh}{{\mathcal K}}

\newcommand{\Mh}{{\mathcal M}}

\newcommand{\Oh}{{\mathcal O}}

\newcommand{\Qh}{{\mathcal Q}}

\newcommand{\Uh}{{\mathcal U}}

\newcommand{\Zh}{{\mathcal Z}}

\newcommand{\ad}{\mathrm{ad}\,}

\newcommand{\au}{\approx_{\mathrm{a.u.}}}

\newcommand{\be}{\mathbf{1}}

\newcommand{\dist}{\mathrm{dist}}

\newcommand{\halb}{\frac{1}{2}}
\newcommand{\her}{\mathrm{her}}

\newcommand{\id}{\mathrm{id}}

\newcommand{\range}{\mathrm{range}\,}

\newcommand{\comp}{\mbox{\scriptsize $\,\circ\,$}}

\newcounter{number}[section]

\newenvironment{nummer}{\refstepcounter{number}{\noindent\arabic{section}.\arabic{number}}}{}

\newcommand{\bn}{\noindent \begin{nummer} \rm}
\newcommand{\en}{\end{nummer}}

\newenvironment{ntheorem}{\noindent {\sc Theorem:} \it}{}
\newenvironment{nlemma}{\noindent {\sc Lemma:} \it}{}
\newenvironment{nprop}{\noindent {\sc Proposition:} \it}{}
\newenvironment{ndefn}{\noindent {\sc Definition:} \it}{}
\newenvironment{ncor}{\noindent {\sc Corollary:} \it}{}

\newenvironment{nremark}{\noindent {\sc Remark:}}{}

\newenvironment{nexamples}{\noindent {\sc Examples:} }{}

\newenvironment{nnotation}{\noindent {\sc Notation:} }{}
\newenvironment{nquestion}{\noindent {\sc Question:} }{}
\newenvironment{nproof}{\noindent {\sc Proof:}}{\mbox{}\hfill 
\rule[-.2ex]{.25em}{1.8ex}}

\title[Strongly self-absorbing $C^{*}$-algebras]{{Strongly self-absorbing $C^{*}$-algebras}}

\dedicatory{Dedicated to George Elliott on the occasion of his 60th birthday.}

\begin{document}

\author{Andrew S.\ Toms}
\address{Department of Mathematics and Statistics, University of New Brunswick\\
Frederic- 
\indent ton, New Brunswick\\
E3B 5A3}

\email{atoms@math.unb.ca}

\author{Wilhelm Winter}
\address{Mathematisches Institut der Universit\"at M\"unster\\
Einsteinstr. 62\\ D-48149 M\"unster}

\email{wwinter@math.uni-muenster.de}

\date{\today}
\thanks{The first author was supported by an NSERC Postdoctoral Fellowship, and the second
author\\
\indent by DFG (through the SFB 478), EU-Network  Quantum Spaces - Noncommutative 
Geometry\\
\indent (Contract No. HPRN-CT-2002-00280).}

\setcounter{section}{-1}

\begin{abstract}

Say that a separable, unital $C^*$-algebra $\Dh \ncong \mathbb{C}$ is 
strongly self-absorbing if there exists an isomorphism $\varphi: \Dh \to \Dh 
\otimes \Dh$ such that $\varphi$ and $\mathrm{id}_{\Dh} \otimes \be_{\Dh}$ are
approximately unitarily
equivalent $*$-homomorphisms.  We study this class of algebras, which includes the Cuntz algebras 
$\mathcal{O}_2$, $\mathcal{O}_{\infty}$, the UHF algebras of infinite type, the Jiang--Su
algebra $\Zh$ and tensor products of $\Oh_{\infty}$ with UHF algebras of infinite type.  Given a strongly self-absorbing $C^{*}$-algebra $\Dh$ we characterise 
when a separable $C^*$-algebra absorbs $\Dh$ tensorially (i.e., is $\Dh$-stable),
and prove closure properties for the class of separable $\Dh$-stable $C^*$-algebras. 
Finally, we compute the possible K-groups and prove a number of classification results which suggest that the examples listed above are the only strongly self-absorbing $C^*$-algebras.   
\end{abstract}

\maketitle

\section{Introduction}\label{intro}

Elliott's program to classify nuclear $C^{*}$-algebras via $\mathrm{K}$-theoretic 
invariants (see \cite{El2} for
an introductory overview) has met with considerable success since his seminal 
classification of approximately
finite-dimensional (AF) algebras via the scaled ordered
$\mathrm{K}_0$-group (\cite{El1}).  An exhaustive list of the contributions to 
this pursuit would be
prohibitively long, but salient works include \cite{El1}, 
\cite{El3}, \cite{EGL}, \cite{JS1}, \cite{K}, \cite{KP}, \cite{Li1}, 
and \cite{M}.  A 
great variety of $C^*$-algebras are studied by these authors, and, 
despite their apparent differences, all of them 
have been classified by $\mathrm{K}$-theoretic invariants.

Upon studying the literature related to Elliott's program, 
one finds that certain $C^*$-algebras have 
been starting points for major stages of the classification program:
UHF algebras in the stably finite case, and the Cuntz algebras in
the purely infinite case.  One can safely say that, among the Cuntz algebras,
$\mathcal{O}_2$ and $\mathcal{O}_{\infty}$ stand out;  they are cornerstones 
of the Kirchberg--Phillips classification of simple purely infinite $C^*$-algebras
and of Kirchberg's classification of non-simple $\Oh_{2}$-absorbing $C^{*}$-algebras 
(in the case where said algebras satisfy the Universal Coefficients Theorem). 
There is evidence that the Jiang--Su algebra $\Zh$, which has recently come to the fore 
of the classification program, plays a role in the stably finite case similar to that of 
$\Oh_{\infty}$ in the purely infinite case (cf.\ \cite{R4}, \cite{TW2} and \cite{Wi5}).

One might reasonably ask whether
there is an abstract property which singles these algebras out from among
their peers.  UHF algebras (at least those of infinite type), $\mathcal{O}_2$,
$\mathcal{O}_{\infty}$ and $\Zh$ are all isomorphic to their tensor squares in
a strong sense;  for each algebra $\Dh$ from this list there exists an isomorphism $\varphi: \Dh \to \Dh 
\otimes \Dh$ such that $\varphi$ and $\mathrm{id}_{\Dh} \otimes \be_{\Dh}$ are
approximately unitarily equivalent $*$-homomorphisms.  In the sequel we refer to
such algebras as strongly self-absorbing whenever they are separable, unital, and
not isomorphic to the complex numbers.
Studying strongly self-absorbing
$C^*$-algebras in the abstract, one finds that the flip automorphisms on their tensor 
squares are approximately inner, whence  they are simple and nuclear by results of \cite{ER}. 
Moreover, they are either purely infinite or stably finite with unique trace (by results 
of Kirchberg, Blackadar and Handelman, and Haagerup). For a strongly self-absorbing $C^*$-algebra
$\Dh$, we say that a second $C^*$-algebra $A$ is $\Dh$-stable if the tensor
product $A \otimes \Dh$ is isomorphic to $A$. Extending results of Kirchberg (cf.\ \cite{K2}), 
we establish permanence properties for the class of $\Dh$-stable $C^*$-algebras
under operations such as taking inductive limits, passing to quotients, 
hereditary subalgebras and ideals, and forming extensions (see also \cite{HW} for results on crossed product $C^{*}$-algebras).

On the other hand, we consider questions relating to the classification program.
We establish classification results for certain strongly self-absorbing $C^*$-algebras;
there is evidence that the examples of said algebras presented in the sequel are
the only such.  A complete classification of the purely infinite strongly self-absorbing 
$C^{*}$-algebras (satisfying the UCT) is given; here it turns out that the only examples 
are $\Oh_{2}$, $\Oh_{\infty}$ and tensor products of $\Oh_{\infty}$ with UHF algebras of 
infinte type. Strongly self-absorbing inductive limits
of recursive subhomogeneous algebras are shown to have the property of slow 
dimension growth in the sense of Phillips (\cite{P3}).  As a corollary, we
show that these are either projectionless or UHF algebras of infinite type. 
Similarly, we conclude that the latter are the only strongly self-absorbing AH algebras and, 
in fact, the only locally type I strongly self-absorbing
$C^*$-algebras of real rank zero.  
In subsequent work we will pay special 
attention to the Jiang--Su algebra $\Zh$, and
the class of $\Zh$-stable $C^*$-algebras (\cite{TW2}).  

We wish to point out that our approach to some extent follows the lines of \cite{ER}, in which 
Effros and  Rosenberg studied $C^{*}$-algebras with approximately inner flip. They derived 
abstract properties (such as nuclearity and simplicity) as well as classification results, namely,  
they showed that the only AF algebras with approximately inner flip are the matroid 
ones (or UHF algebras, in the unital case). At that time these were the only known examples 
of such $C^{*}$-algebras. 

\noindent
\emph{Acknowledgements.}
We would like to thank George Elliott for suggesting the study of general self-absorbing
$C^{*}$-algebras, and Eberhard Kirchberg and Mikael R{\o}rdam for many inspiring discussions.  We are also indebted to the referee for a number of helpful comments.

\section{Strongly self-absorbing $C^{*}$-algebras}\label{ssa}

A $C^{*}$-algebra is usually referred to as being self-absorbing if it is isomorphic to its 
tensor product with itself. In general this statement requires specification of a particular 
tensor product. Since we are mostly interested in the nuclear case, there will be no loss in 
generality if we consider only minimal $C^{*}$-algebraic tensor products. 

Self-absorbing $C^{*}$-algebras can be easily constructed: If $A$ is any nuclear and unital 
$C^{*}$-algebra, let $A^{\otimes \infty}$ denote the $C^{*}$-limit of the inductive system 
\[
(A^{\otimes n}, \id_{A^{\otimes n}} \otimes \be_{A})_{n \in \N} \, .
\]
It is not hard to see that $A^{\otimes \infty}$ is self-absorbing. Repeating this process 
will not yield anything new; we have $(A^{\otimes \infty})^{\otimes \infty} \cong A^{\otimes \infty}$.

The examples that motivated this article (see \ref{absorbing-examples}) are self-absorbing 
in a much stronger sense. In this section, we describe the concept of being strongly self-absorbing 
and a number of characterisations and structural properties. First, we recall the notion of and 
some well-known facts about approximate unitary equivalence.

\bn
\begin{ndefn}
For $i=1,2$, let $\varphi_{i}: A \to B$ be a c.p.c.\ map between separable $C^{*}$-algebras. 
We say $\varphi_{1}$ and $\varphi_{2}$ are approximately unitarily (a.u.) equivalent, 
$\varphi_{1} \au \varphi_{2}$, if there is a sequence $(v_{n})_{\N}$ of unitaries in 
the multiplier algebra $\Mh(B)$ such that $\| v_{n}^{*} \varphi_{1}(a)v_{n} - \varphi_{2}(a)\| \stackrel{n \to \infty}{\longrightarrow} 0 \; \forall \, a \in A$.
\end{ndefn}
\en

\bn
\label{transitivity} 
\begin{nprop}
Let $A$, $B$, $C$ and $D$ be separable $C^{*}$-algebras, $C$ and $D$ unital. Suppose $\varphi:A \to B$,  $\alpha,\beta, \gamma:B \to C$,  and $\psi:C \to D$ are $*$-homomorphisms, $\psi$ unital.
\begin{itemize}
\item[(i)] If $\alpha \au \beta$ and $\beta \au \gamma$, then $\alpha \au \gamma$. In other words,  approximate unitary equivalence is a transitive relation.
\item[(ii)] If $\alpha \au \beta$, then $\psi \circ \alpha \au \psi \circ \beta$ and $\alpha \circ \varphi \au \beta \circ \varphi$.
\item[(iii)] Suppose $\alpha$ and $\beta$ are pointwise limits of 
sequences of $*$-homomorphisms $\alpha_{n}, \beta_{n}: B \to C$, respectively. 
If $\alpha_{n} \au \beta_{n}$ for each $n \in \N$, then $\alpha \au \beta$.
\end{itemize}
\end{nprop}
\en

\bn
\label{d-half-flip}
\begin{ndefn}
Let $\Dh$ be a separable unital $C^{*}$-algebra.
\begin{itemize}
\item[(i)] By the flip on the minimal tensor product $\Dh \otimes \Dh$ we mean the automorphism $\sigma_{\Dh}$ of $\Dh \otimes \Dh$ given by $\sigma_{\Dh}(a \otimes b):= b \otimes a, \, a, b \in \Dh$.
\item[(ii)] $\Dh$ is said to have approximately inner flip, if $\sigma_{\Dh}$ is approximately 
unitarily equivalent to the identity map, i.e., $\sigma_{\Dh} \au \id_{\Dh \otimes \Dh}$.
\item[(iii)] $\Dh$ is said to have approximately inner half flip, if the two natural inclusions  
of $\Dh$ into $\Dh \otimes \Dh$ as the first and second factor, respectively, are approximately unitarily equivalent, i.e., $\id_{\Dh} \otimes \be_{\Dh} \au \be_{\Dh} \otimes \id_{\Dh}$.
\item[(iv)] $\Dh$ is strongly self-absorbing, if $\Dh \ncong \C$ and  there is an isomorphism $\varphi: \Dh \to \Dh \otimes \Dh$ satisfying $\varphi \au \id_{\Dh} \otimes \be_{\Dh}$.
\end{itemize}
\end{ndefn}
\en

\bn
\label{r-half-flip}
\begin{nremark}
In Definition \ref{d-half-flip}(iv) we could as well have asked for an isomorphism 
$\psi: \Dh \to \Dh \otimes \Dh$ satisfying $\psi \au  \be_{\Dh} \otimes \id_{\Dh}$. 
Both definitions are equivalent, since, given $\varphi$, we may choose $\psi:= 
\sigma_{\Dh} \circ \varphi$ and vice versa.  
\end{nremark}
\en

\bn
\label{self-absorbing-half-flip}
The preceding remark shows that Definition \ref{d-half-flip}(iv) is in fact symmetric, 
although it is not stated this way. Even more is true:  In Corollary \ref{self-absorbing-flip} 
it will turn out  that strongly self-absorbing $C^{*}$-algebras have approximately inner flip. 
As a first step, we show that they have approximately inner half flip, from which already 
follows that $\varphi$ and $\sigma_{\Dh} \circ \varphi$ are approximately unitarily equivalent. 

\begin{nprop}
If $\Dh$ is separable, unital and strongly self-absorbing, then $\Dh$ has approximately 
inner half flip.
\end{nprop}

\begin{nproof}
Suppose that $\varphi: \Dh \to \Dh \otimes \Dh$ is an isomorphism such that 
\[
\varphi \au \id_{\Dh} \otimes \be_{\Dh} \, .
\]
Define a unital $*$-homomorphism $\psi: \Dh \to \Dh$ by 
\[
\psi:= \varphi^{-1} \circ (\be_{\Dh} \otimes \id_{\Dh}) \, .
\]
Note that
\[
\begin{array}{lll}
\be_{\Dh} \otimes \id_{\Dh} & = & \varphi \circ \varphi^{-1} \circ (\be_{\Dh} \otimes \id_{\Dh}) \\
& = & \varphi \circ \psi \\
& \stackrel{\ref{transitivity}(ii)}{\au} & (\id_{\Dh} \otimes \be_{\Dh}) \circ \psi \\
& = & \psi \otimes \be_{\Dh} \, ;
\end{array}
\]
from Proposition \ref{transitivity}(ii) we also see that 
\[
\begin{array}{lll}
\id_{\Dh} \otimes \be_{\Dh} & = & \sigma_{\Dh} \circ (\be_{\Dh} \otimes \id_{\Dh}) \\
& \au & \sigma_{\Dh} \circ (\psi \otimes \be_{\Dh}) \\
& = & \be_{\Dh} \otimes \psi \, .
\end{array}
\]
We now proceed to obtain
\[
\begin{array}{lll}
\psi \otimes \be_{\Dh} & = & (\id_{\Dh} \otimes \varphi^{-1}) \circ (\psi \otimes \be_{\Dh} \otimes \be_{\Dh}) \\
& \au & (\id_{\Dh} \otimes \varphi^{-1}) \circ (\be_{\Dh} \otimes \id_{\Dh} \otimes \be_{\Dh}) \\
& \au & (\id_{\Dh} \otimes \varphi^{-1}) \circ (\be_{\Dh} \otimes \be_{\Dh} \otimes \psi) \\
& \au & (\id_{\Dh} \otimes \varphi^{-1}) \circ (\id_{\Dh} \otimes \be_{\Dh} \otimes \be_{\Dh})\\
& = & \id_{\Dh} \otimes \be_{\Dh} \, .
\end{array}
\]
By transitivity of approximate unitary equivalence we thus have 
\[
\be_{\Dh} \otimes \id_{\Dh} \au \id_{\Dh} \otimes \be_{\Dh} \, .
\]
\end{nproof}
\en

\bn
\label{simple-nuclear-half-flip}
Before continuing our analysis of strongly self-absorbing $C^{*}$-algebras, we recall an 
important structure result about $C^{*}$-algebras with approximately inner half flips. 
In the case where $\Dh$ has approximately inner flip the statement was already observed 
in \cite{ER}. In the form we state below, the result was shown in \cite{KP}. 

\begin{ntheorem}
If a separable unital $C^{*}$-algebra $\Dh$ has approximately inner half flip, it is simple and nuclear.
\end{ntheorem}
\en

\bn
\label{dichotomy}
The following result provides a first instance why strongly self-absorbing $C^{*}$-algebras 
play an important role in the classification program.

\begin{ntheorem}
A separable unital strongly self-absorbing $C^{*}$-algebra $\Dh$ is either purely infinite 
or stably finite with a unique tracial state.
\end{ntheorem}

\begin{nproof}
The fact that $\Dh \cong \Dh \otimes \Dh$ is either stably finite or purely infinite 
is due to Kirchberg (cf.\ \cite{R1}, Theorem 4.1.10). That $\Dh$ admits a tracial 
state follows from results of Blackadar, Handelman and of Haagerup (cf.\ \cite{R1}, Theorem 1.1.4). 
The tracial state has to be unique by (literally the same proof as that of) \cite{ER}, 
Proposition 2.10; the argument applies since $\Dh$ has approximately inner half flip.
\end{nproof}
\en

\bn
\label{flip-products}
\begin{nprop}
If $\Dh$ and $\Eh$ are separable unital $C^{*}$-algebras both with approximately inner 
(half) flips, then $\Dh \otimes \Eh$ also has approximately inner (half) flip. If $\Dh$ and $\Eh$ are strongly self-absorbing, then so is $\Dh \otimes \Eh$.
\end{nprop}

\begin{nproof}
When $\Dh$ and $\Eh$ have approximately inner flips, this was shown in \cite{ER}, 
Corollary 2.4 (to Lemma 2.3). In the case of approximately inner half flips, the 
same proof applies almost verbatim, only replacing the identy maps and the flip 
automorphisms on $\Dh \otimes \Dh$ and $\Eh \otimes \Eh$ by the appropriate canonical 
unital embeddings of $\Dh$ and $\Eh$ into $\Dh \otimes \Dh$ and $\Eh \otimes \Eh$, respectively. The last statement is proved similarly.
\end{nproof}
\en

\bn
\label{infinite-products}
\begin{nprop}
Suppose $\Dh$ is separable and unital and $\Dh$ has approximately inner half flip. Then: 
\begin{itemize}
\item[(i)] $\Dh^{\otimes \infty}$ has approximately inner flip.
\item[(ii)] $\Dh^{\otimes \infty}$ is strongly self-absorbing.
\item[(iii)] There is a sequence of $*$-homomorphisms 
\[
(\varphi_{n}: \Dh^{\otimes \infty} \otimes \Dh^{\otimes \infty} \to \Dh^{\otimes \infty})_{n \in \N} 
\]  
satisfying 
\[
\|\varphi_{n}(d \otimes \be_{\Dh^{\otimes \infty}}) - d\| \stackrel{n \to \infty}{\longrightarrow} 0 \; \forall \, d \in \Dh^{\otimes \infty} \, .
\]
\end{itemize}
\end{nprop}

\begin{nproof}
(i) By the definition of $\Dh^{\otimes \infty}$ as an inductive limit it clearly suffices to show that, for $k \in \N$, we have
\[
\lambda_{k} \au \lambda_{k} \circ \sigma_{\Dh^{\otimes k}} \, ,
\]
where $\lambda_{k}: \Dh^{\otimes k} \otimes \Dh^{\otimes k} \to \Dh^{\otimes 2k} \otimes \Dh^{\otimes 2k}$ is given by
\[
\lambda_{k} = (\id_{\Dh^{\otimes k}} \otimes \be_{\Dh^{\otimes k}}) \otimes  (\id_{\Dh^{\otimes k}} \otimes \be_{\Dh^{\otimes k}}) 
\]
and $\sigma_{\Dh^{\otimes k}}$ is the flip on $\Dh^{\otimes k} \otimes \Dh^{\otimes k}$. 

We denote the embedding of $\Dh^{\otimes k}$ into $(\Dh^{\otimes k})^{\otimes 4}$ as the $i$-th factor by $\iota_{k}^{(i)}$. Then, we define $*$-homomorphisms
\[
\iota_{k}^{(i,j)}: (\Dh^{\otimes k})^{\otimes 2} \to (\Dh^{\otimes k})^{\otimes 4}\, , \, i\neq j \in \{1, \ldots, 4\} \, ,
\]
by 
\[
\iota^{(i,j)}_{k}|_{\Dh^{\otimes k} \otimes \be_{\Dh^{\otimes k}}} = \iota^{(i)} \mbox{ and } 
\iota^{(i,j)}_{k}|_{\be_{\Dh^{\otimes k}} \otimes \Dh^{\otimes k}} = \iota^{(j)} \, .
\]
Note that $\iota_{k}^{(i,j)}$ is well-defined this way, since $\iota^{(i)}(\Dh^{\otimes k})$ and $\iota^{(j)}(\Dh^{\otimes k})$ commute.
Identifying $\Dh^{\otimes 2k} \otimes \Dh^{\otimes 2k}$ with $(\Dh^{\otimes k})^{\otimes 4}$ in the obvious way, with these definitions we have
\[
\lambda_{k} = \iota_{k}^{(1,3)} \mbox{ and } \lambda_{k} \circ \sigma_{\Dh^{\otimes k}} = \iota_{k}^{(3,1)} \, .
\]

Now let $i,j,l \in \{1, \ldots,4\}$ be pairwise distinct. By Proposition \ref{flip-products}, 
$\Dh^{\otimes k}$ has approximately inner half flip, so there is a sequence $(v_{m})_{m \in \N}$ 
of unitaries in $\Dh^{\otimes k} \otimes \Dh^{\otimes k}$ intertwining the two canonical embeddings 
of $\Dh^{\otimes k}$ into $\Dh^{\otimes k} \otimes \Dh^{\otimes k}$. But then 
$(\iota^{(j,l)}_{k}(v_{m}))_{m \in \N} \subset (\Dh^{\otimes k})^{\otimes 4}$ is 
a sequence of unitaries intertwining $\iota^{(j)}_{k}$ and $\iota^{(l)}_{k}$; 
the $\iota^{(j,l)}_{k}(v_{m})$ commute with $\iota^{(i)}_{k}(\Dh^{\otimes k})$, 
whence they even intertwine $\iota^{(i,j)}_{k}$ and $\iota^{(i,l)}_{k}$.  Therefore we have 
\[
\iota_{k}^{(i,j)} \au \iota_{k}^{(i,l)}
\]
and, similarly,
\[
\iota_{k}^{(i,j)} \au \iota_{k}^{(l,j)} \, .
\]
In particular we obtain 
\[
\iota_{k}^{(1,3)} \au \iota_{k}^{(1,2)} \au \iota_{k}^{(3,2)} \au \iota_{k}^{(3,1)}
\]
and, by transitivity of a.u.\ equivalence,
\[
\lambda_{k} = \iota_{k}^{(1,3)} \au \iota_{k}^{(3,1)} = \lambda_{k} \circ \sigma_{\Dh^{\otimes k}} \, .
\]
(ii) For $k \in \N$ define $*$-homomorphisms $\alpha_{k}:\Dh^{\otimes k} \to \Dh^{\otimes k+1}$ by 
\[
\alpha_{k}= \id_{\Dh^{\otimes k}} \otimes \be_{\Dh} \, .
\]
Then we have 
\[
\Dh^{\otimes \infty} = \lim_{\to} (\Dh^{\otimes k},\alpha_{k}) \, ,
\]
but we also obtain
\begin{equation}
\label{embedding1}
\Dh^{\otimes \infty} \otimes \Dh^{\otimes \infty} = \lim_{\to} (\Dh^{\otimes k} \otimes \Dh^{\otimes k},\alpha_{k} \otimes \alpha_{k})
\end{equation}
and
\begin{equation}
\label{embedding2}
\Dh^{\otimes \infty} = \lim_{\to} (\Dh^{\otimes 2k} ,\alpha_{2k+1} \circ \alpha_{2k}) \, .
\end{equation}
Next, define isomorphisms 
\[
\psi_{k}: \Dh^{\otimes 2k} \to \Dh^{\otimes k} \otimes \Dh^{\otimes k}
\]
by
\[
\psi_{k}(a_{1} \otimes b_{1} \otimes \ldots \otimes a_{k} \otimes b_{k}) = a_{1} \otimes \ldots \otimes a_{k} \otimes b_{1} \otimes \ldots \otimes b_{k} \, .
\]
We have 
\[
 \psi_{k+1} \circ \alpha_{2k+1} \circ \alpha_{2k} = (\alpha_{k} \otimes \alpha_{k}) \circ \psi_{k} \, ,
\]
so, by (\ref{embedding1}), (\ref{embedding2}) and the universal property of inductive limits, we see that the $\psi_{k}$ induce a $*$-homomorphism 
\[
\psi: \Dh^{\otimes \infty} \to \Dh^{\otimes \infty} \otimes \Dh^{\otimes \infty} \, .
\]
Since each $\psi_{k}$ is an isomorphism, so is $\psi$.

We want to prove that $\psi \au \id_{\Dh^{\otimes \infty}} \otimes \be_{\Dh^{\otimes \infty}}$. Since 
\[
\Dh^{\otimes \infty} \otimes \Dh^{\otimes \infty} = \lim_{\to} (\Dh^{\otimes 2^{m}} \otimes \Dh^{\otimes 2^{m}}, \lambda_{2^{m}}) 
\]
and
\[
\Dh^{\otimes \infty} = \lim_{\to} (\Dh^{\otimes 2^{m}}, \id_{\Dh^{\otimes 2^{m}}} \otimes \be_{\Dh^{\otimes 2^{m}}}) \, ,
\]
it will be enough to show that
\[
\lambda_{k} \circ \psi_{k} \circ (\id_{\Dh^{\otimes k}} \otimes \be_{\Dh^{\otimes k}}) \au \lambda_{k} \circ (\id_{\Dh^{\otimes k}} \otimes \be_{\Dh^{\otimes k}})
\]
for any $k \in \N$ (in fact, it would suffice to show the above for each $k \in \{2^{m}\, | \, m \in \N\}$). First, define $*$-homomorphisms 
\[
\beta_{k}: (\Dh^{\otimes k})^{4} \to (\Dh^{\otimes k})^{4}\]
by
\[
\beta_{k}:= \id_{\Dh^{\otimes k}} \otimes \sigma_{\Dh^{\otimes k}} \otimes \id_{\Dh^{\otimes k}} \, .
\]
Then, again identifying $\Dh^{\otimes 2k} \otimes \Dh^{\otimes 2k}$ with $(\Dh^{\otimes k})^{\otimes 4}$, for $i=1,2$ we obtain $*$-homomorphisms
\[
\beta_{k} \circ (\psi_{k} \otimes \id_{\Dh^{\otimes 2k}}) \circ \iota_{k}^{(i)}: \Dh^{\otimes k} \to \Dh^{\otimes 2k} \otimes \Dh^{\otimes 2k} \, ;
\]
one easily checks that
\[
\beta_{k} \circ (\psi_{k} \otimes \id_{\Dh^{\otimes 2k}}) \circ \iota_{k}^{(1)} = \lambda_{k} \circ \psi_{k} \circ (\id_{\Dh^{\otimes k}} \otimes \be_{\Dh^{\otimes k}})
\]
and
\[
\beta_{k} \circ (\psi_{k} \otimes \id_{\Dh^{\otimes 2k}}) \circ \iota_{k}^{(3)} = \iota^{(2)}_{k} \, .
\]
Since $\iota_{k}^{(i)} \au \iota_{k}^{(j)}$, by  Proposition \ref{transitivity}(ii) we have
\[
\begin{array}{rll}
\lambda_{k} \circ \psi_{k} \circ (\id_{\Dh^{\otimes k}} \otimes \be_{\Dh^{\otimes k}}) &= &\beta_{k} \circ (\psi_{k} \otimes \id_{\Dh^{\otimes 2k}}) \circ \iota_{k}^{(1)} \\
& \au & \beta_{k} \circ (\psi_{k} \otimes \id_{\Dh^{\otimes 2k}}) \circ \iota_{k}^{(3)} \\
&= &\iota_{k}^{(2)} \\
& \au & \iota_{k}^{(1)} \\
& = & \lambda_{k} \circ (\id_{\Dh^{\otimes k}} \otimes \be_{\Dh^{\otimes k}}) \, .
\end{array}
\]
(iii) If $(u_{n})_{\N} \subset \Dh^{\otimes \infty}$ is a sequence of unitaries such that 
\[
\|u_{n}^{*}\psi(d)u_{n} - d \otimes \be_{\Dh^{\otimes \infty}}\| \stackrel{n \to \infty}{\longrightarrow} 0 \; \forall \, d \in \Dh^{\otimes \infty}\, ,
\]
then
\[
\| \psi^{-1} (u_{n} (d \otimes \be_{\Dh^{\otimes \infty}}) u_{n}^{*}) - d\| \stackrel{n \to \infty}{\longrightarrow}  0 \; \forall \, d \in \Dh^{\otimes \infty} \, .
\]
Therefore, we may define the $*$-homomorphisms $\varphi_{n}$ by
\[
\varphi_{n}(d_{1} \otimes d_{2}):= \psi^{-1} (u_{n}(d_{1} \otimes d_{2}) u_{n}^{*}) \,, \; d_{1}, d_{2} \in \Dh^{\otimes \infty}.
\]
\end{nproof}
\en

\bn
\label{t-self-absorbing}
\begin{nprop}
Let $\Dh$ be a separable unital $C^{*}$-algebra such that $\Dh$ has approximately inner half flip. 
Then $\Dh$ is strongly self-absorbing iff one of the following equivalent conditions holds:
\begin{itemize}
\item[(i)] There exists a unital $*$-homomorphism $\gamma:\Dh \otimes \Dh \to \Dh$ which satisfies $\gamma \circ (\id_{\Dh} \otimes \be_{\Dh}) \au \id_{\Dh}$. 
\item[(ii)] There are a unital $*$-homomorphism $\gamma:\Dh \otimes \Dh \to \Dh$ and an approximately central sequence of unital endomorphisms of $\Dh$.
\item[(iii)] There exists an approximately central sequence of unital $*$-homomorphisms $\Dh^{\otimes \infty} \to \Dh$.
\item[(iv)] There exists an isomorphism $\Dh \to \Dh^{\otimes \infty}$.
\end{itemize}
\end{nprop}

\begin{nproof}
If $\Dh$ is strongly self-absorbing, there are an isomorphism $\varphi: \Dh \to \Dh \otimes \Dh$ and a sequence of unitaries $(u_{n})_{\N} \subset \Dh \otimes \Dh$ such that 
\[
\|u_{n}^{*}\varphi(d)u_{n} - d \otimes \be_{\Dh}\| \stackrel{n \to \infty}{\longrightarrow} 0 \; \forall \, d \in \Dh \, .
\]
Set $\gamma:= \varphi^{-1}$, then $(\gamma(u_{n}))_{\N} \subset \Dh$ is a sequence of unitaries satisfying
\[
\|\gamma(u_{n}^{*})d \gamma(u_{n}) - \gamma(d \otimes \be_{\Dh})\| \stackrel{n \to \infty}{\longrightarrow} 0 \; \forall \, d \in \Dh \, ,
\]
so (i) holds.

\noindent
(i) $\Rightarrow$ (ii): Let $(v_{n})_{\N} \subset \Dh$ be a sequence of unitaries such that 
\[
\|v_{n}^{*} \gamma(d \otimes \be_{\Dh}) v_{n} - d \| \stackrel{n \to \infty}{\longrightarrow} 0 \; \forall \, d \in \Dh \, .
\]
Define $*$-homomorphisms $\varphi_{n}: \Dh \to \Dh$ by 
\[
\varphi_{n}(d):= v_{n}^{*} \gamma(\be_{\Dh} \otimes d)v_{n} \, ;
\]
it then follows from $[\varphi_{n}(d_{1}), v_{n}^{*}\gamma(d_{2} \otimes \be_{\Dh})v_{n}] = 0$ that 
\[
\|[\varphi_{n}(d_{1}),d_{2}] \stackrel{n \to \infty}{\longrightarrow} 0 \; \forall \, d_{1}, d_{2} \in \Dh \, .
\]

\noindent
(ii) $\Rightarrow$ (iii): It obviously suffices to construct a unital $*$-homomorphism 
\[
\psi: \Dh^{\otimes \infty} \to \Dh \, .
\]
For $k \in \N$, define unital $*$-homomorphisms $\gamma_{k}: \Dh^{\otimes k+1} \to \Dh^{\otimes k}$ by
\[
\gamma_{k}:= \id_{\Dh^{\otimes k-1}} \otimes \gamma
\]
and $\psi_{k}: \Dh^{\otimes k} \to \Dh$ by 
\[
\psi_{k}:= \gamma_{1} \circ \cdots \circ \gamma_{k+1} \circ (\id_{\Dh^{\otimes k}} \otimes \be_{\Dh^{\otimes 2}}) \, .
\]
We now have
\begin{equation*}
\psi_{k}  =  \gamma_{1} \circ \ldots \circ \gamma_{k+1} \circ \gamma_{k+2} \circ (\id_{\Dh^{\otimes k}} \otimes \be_{\Dh^{\otimes 3}}) =  \psi_{k+1} \circ (\id_{\Dh^{\otimes k}} \otimes \be_{\Dh}) \, ,
\end{equation*}
from which follows that the $\psi_{k}$ induce a (unital) $*$-homomorphism $\psi: \Dh^{\otimes \infty} \to \Dh$.

\noindent
(iii) $\Rightarrow$ (iv): By Proposition \ref{infinite-products}, $\Dh^{\otimes \infty}$ 
is strongly self-absorbing, and it follows from \cite{R1}, Theorem 7.2.2 (cf.\ also Theorem 
\ref{rordams-intertwining} below) that $\Dh \cong \Dh \otimes \Dh^{\otimes \infty}$. We 
obtain an isomorphism $\Dh \otimes \Dh^{\otimes \infty} \cong \Dh^{\otimes \infty}$ from 
the right shift of the inductive system defining $\Dh^{\otimes \infty}$.

Finally, if $\Dh \cong \Dh^{\otimes \infty}$, then $\Dh$ is strongly self-absorbing by Proposition \ref{infinite-products}.
\end{nproof}
\en

\bn
\label{self-absorbing-flip}
\begin{ncor}
If $\Dh$ is separable, unital and strongly self-absorbing, then $\Dh \cong \Dh^{\otimes k} 
\cong \Dh^{\otimes \infty}$ for any $k \in \N$ and $\Dh$ has approximately inner flip.  
\end{ncor}

\begin{nproof}
That $\Dh \cong \Dh^{\otimes k}$ for any $k$ is trivial; the other statements simply summarize \ref{infinite-products}(i) and \ref{t-self-absorbing}(iv).
\end{nproof}
\en

\bn
\label{equivalent-endomorphisms}
\begin{ncor}
Let $A$ and $\Dh$ be separable unital $C^{*}$-algebras, with $\Dh$ strongly self-absorbing. 
Then, any two unital  $*$-homomorphisms $\alpha,\beta:\Dh \to A \otimes \Dh$ are a.u.\ equivalent. 
In particular, any two unital endomorphisms of $\Dh$ are a.u.\ equivalent.  
\end{ncor}

\begin{nproof}
By Proposition \ref{infinite-products}(iii) (in connection with Proposition \ref{t-self-absorbing}(iv)) 
there is a sequence of unital $*$-homomorphisms $\varphi_{n}: \Dh \otimes \Dh \to \Dh$ such that 
$\varphi_{n} \comp (\id_{\Dh} \otimes \be_{\Dh}) \to \id_{\Dh}$ pointwise.

For each $n \in \N$ we define unital $*$-homomorphisms $\bar{\alpha}_{n}, \bar{\beta}_{n}: \Dh \to A \otimes \Dh$ by
\[
\bar{\alpha}_{n}:= (\id_{A} \otimes \varphi_{n}) \circ (\alpha \otimes \id_{\Dh}) \circ (\id_{\Dh} \otimes \be_{\Dh})
\]
and 
\[
\bar{\beta}_{n}:= (\id_{A} \otimes \varphi_{n}) \circ (\beta \otimes \id_{\Dh}) \circ (\id_{\Dh} \otimes \be_{\Dh}) \, .
\]
From Proposition \ref{transitivity}(ii) we see that 
\[
\begin{array}{cll}
\bar{\alpha}_{n} & \au &  (\id_{A} \otimes \varphi_{n}) \circ (\alpha \otimes \id_{\Dh}) \circ (\be_{\Dh} \otimes \id_{\Dh})\\
& = & (\id_{A} \otimes \varphi_{n}) \circ (\be_{A} \otimes \be_{\Dh} \otimes \id_{\Dh}) \\
& = & (\id_{A} \otimes \varphi_{n}) \circ (\beta \otimes \id_{\Dh}) \circ (\be_{\Dh} \otimes \id_{\Dh})\\
& \au & \bar{\beta}_{n} 
\end{array}
\]
and, from \ref{transitivity}(i), we obtain $\bar{\alpha}_{n} \au \bar{\beta}_{n} \; \forall \, n \in \N$. 
But we obviously have $\bar{\alpha}_{n} \to \alpha$ and $\bar{\beta}_{n} \to \beta$ pointwise, so $\alpha \au \beta$ by Proposition \ref{transitivity}(iii). 
The second statement follows with $A = \Dh$, since $\Dh \cong \Dh \otimes \Dh$ by assumption.
\end{nproof}
\en

\bn
\label{unitary-homotopies}
Recall that a unital $C^{*}$-algebra is said to be $\mathrm{K}_{1}$-injective, if the canonical homomorphism $\Uh(\Dh)/\Uh_{0}(\Dh) \to \mathrm{K}_{1}(\Dh)$ is injective.

\begin{nprop}
Let $\Dh$ be a separable, unital, strongly self-absorbing $C^{*}$-algebra. Then the unitaries 
implementing the approximately inner flip on $\Dh \otimes \Dh$ may be chosen to represent $0$ 
in $\mathrm{K}_{1}(\Dh \otimes \Dh) \cong \mathrm{K}_{1}(\Dh)$. In particular, if $\Dh$ is 
$\mathrm{K}_{1}$-injective, then the unitaries may be chosen to be homotopic to $\be_{\Dh}$.
\end{nprop}

\begin{nproof}
For $k \in \N$ define $*$-homomorphisms $\lambda_{k}: \Dh^{\otimes 2^{k}} \otimes \Dh^{\otimes 2^{k}} \to \Dh^{\otimes 2^{k+1}} \otimes \Dh^{\otimes 2^{k+1}}$ by
\[
\lambda_{k} = (\id_{\Dh^{\otimes 2^{k}}} \otimes \be_{\Dh^{\otimes 2^{k}}}) \otimes  (\id_{\Dh^{\otimes 2^{k}}} \otimes \be_{\Dh^{\otimes 2^{k}}}) \, .
\]
We then have 
\[
\Dh^{\otimes \infty} \otimes \Dh^{\otimes \infty} = \lim_{\to} (\Dh^{\otimes 2^{k}} \otimes \Dh^{\otimes 2^{k}}, \lambda_{k}) \, ;
\]
denote the canonical embedding of $\Dh^{\otimes 2^{k}} \otimes \Dh^{\otimes 2^{k}}$ 
into $\Dh^{\otimes \infty} \otimes \Dh^{\otimes \infty}$ by $\lambda_{k,\infty}$. 

By Corollary \ref{self-absorbing-flip}, $\Dh$ has approximately inner flip, and so 
has $\Dh^{\otimes 2^{k}}$ by  Proposition \ref{flip-products}.  Let  $(u_{k,n})_{n \in \N} 
\subset \Dh^{\otimes 2^{k}} \otimes \Dh^{\otimes 2^{k}}$ be a sequence of unitaries approximating 
the flip $\sigma_{\Dh^{\otimes 2^{k}}}$ on  $\Dh^{\otimes 2^{k}} \otimes \Dh^{\otimes 2^{k}}$. 
Choosing a suitable sequence $(n_{k})_{k \in \N} \subset \N$,  it is then not hard to  
obtain unitaries $u_{k,n_{k}} \in \Dh^{\otimes 2^{k}} \otimes \Dh^{\otimes 2^{k}}$,  such that for any $m \in \N$ we have
\[
\lambda_{k,\infty}(u_{k,n_{k}}^{*})  \lambda_{m, \infty}(c \otimes d) \lambda_{k,\infty}(u_{k,n_{k}}) 
\stackrel{k \to \infty}{\longrightarrow}  \lambda_{m, \infty}(d \otimes c) \; \forall \, c,d \in \Dh^{\otimes 2^{m}} \, .
\]
But then it is obvious that for any $m \in \N$ we also have 
\begin{equation*}
\lambda_{k+1,\infty}(u_{k,n_{k}}^{*} \otimes u_{k,n_{k}})  \lambda_{m, \infty}(c \otimes d) 
\lambda_{k+1,\infty}(u_{k,n_{k}}\otimes u_{k,n_{k}}^{*}) \stackrel{k \to \infty}{\longrightarrow} \lambda_{m, \infty}(d \otimes c) 
\end{equation*}
for all $c,d \in \Dh^{\otimes \infty}$, so we may define elements
\[
v_{k}:= \lambda_{k+1,\infty}(u_{k,n_{k}}\otimes u_{k,n_{k}}^{*}) \, ,
\]
which form a sequence of unitaries in $\Dh^{\otimes \infty} \otimes \Dh^{\otimes \infty}$ 
approximating the flip $\sigma_{\Dh^{\otimes \infty}}$. 

Let $u \in \Dh$ be any unitary element.  Then we have 
\[
v_k(u \otimes \be_{\Dh})v_k^*(\be_{\Dh} \otimes u^*) \stackrel{k \to \infty}{\longrightarrow} \be_{\Dh \otimes \Dh} \equiv \be_{\Dh} \, .
\]
On the other hand, 
\[
\mathrm{K}_1(v_k(u \otimes \be_{\Dh})v_k^*(\be_{\Dh} \otimes u^*)) = \mathrm{K}_1(u \otimes u^*) \; \forall \, k \, ,
\]
whence $\mathrm{K}_{1}(u \otimes u^{*}) = 0$. Since $u$ was arbitrary, we obtain $\mathrm{K}_{1}(u_{k,n_{k}} \otimes u_{k, n_{k}}^{*}) = 0$ for any (fixed) $k$; we conclude that $\mathrm{K}_1(v_k) = 0$, $\forall \, k \in \mathbb{N}$, as desired.
\end{nproof}
\en

\bn
\label{absorbing-examples}
\begin{nexamples}
(i) Recall that a UHF algebra is usually written in the form 
\[
M_{p_{1}}^{\otimes k_{1}} \otimes M_{p_{2}}^{\otimes k_{2}} \otimes \ldots
\]
or, more convenient, as $M_{\mathfrak n}$, for a supernatural number ${\mathfrak n}= 
p_{1}^{k_{1}} \cdot p_{2}^{k_{2}} \ldots$, where $(p_{i})_{i \in \N}$ is an enumeration of 
the primes and the exponents $k_{i}$ are in $\N \cup \{\infty\}$ ($\N$ containing zero).

It follows from elementary linear algebra that the flip automorphism on $M_{r} \otimes M_{r}$ 
is inner for any $r \in \N$. As a consequence, any UHF algebra has approximately inner flip. 
Proposition \ref{infinite-products} and Corollary \ref{self-absorbing-flip} now show that $B=  M_{\mathfrak n}$ 
is strongly self-absorbing iff $k_{i} \in \{0, \infty\} \; \forall \, i \in \N$ (and at least 
one $k_{i}$ is nonzero). In other words, ${\mathfrak n}$ is nontrivial and each prime that 
occurs in ${\mathfrak n}$ has to occur infinitely many times. We will call such a $B$ a UHF 
algebra of infinite type. It is obvious that $B$ is of infinite type if and only if it is self-absorbing in the ordinary sense.

In \cite{ER} it was shown that the only unital AF algebras with approximately inner flip are UHF algebras. 
Below we will prove a similar result, namely that the only unital strongly self-absorbing AH algebras are UHF of infinite type. 

Recall also that it is well-known that $M_{r}$, and hence any UHF algebra, is $\mathrm{K}_{1}$-injective.

\noindent 
(ii) In \cite{Cu1}, Cuntz introduced the algebras $\Oh_{n}$ (with $n\in \{2,3, \ldots\} \cup 
\{\infty\}$). These are universal $C^{*}$-algebras generated by $n$ isometries with certain 
relations; they are nuclear, simple and  purely infinite. In Kirchberg's  classification of 
purely infinite simple nuclear $C^{*}$-algebras, the Cuntz algebras $\Oh_{2}$ and $\Oh_{\infty}$ 
play a particularly important role. By results of Elliott and of R{\o}rdam (cf.\ \cite{R1}), 
they have approximately inner half flip and satisfy $\Oh_{2} \cong \Oh_{2}^{\otimes \infty} $ 
and $\Oh_{\infty} \cong \Oh_{\infty}^{\otimes \infty}$, respectively, so they are strongly 
self-absorbing by Proposition \ref{t-self-absorbing}. It follows from \cite{Cu2} that $\Oh_{2}$ 
and $\Oh_{\infty}$ are both $\mathrm{K}_{1}$-injective.

\noindent
(iii) If $B$ is a UHF algebra of infinite type, then $B \otimes \Oh_{\infty}$ is strongly 
self-absorbing by Proposition \ref{flip-products}. Since $B \otimes \Oh_{\infty}$ is simple and purely 
infinite, it is $\mathrm{K}_{1}$-injective by \cite{Cu2}. 

Note that $B$ and $B \otimes \Oh_{\infty}$ are KK-eqivalent, but are not isomorphic ($B$ 
is stably finite, whereas $B \otimes \Oh_{\infty}$ is purely infinite). If $B_{1}$ and 
$B_{2}$ are UHF algebras of infinite type, then by Kirchberg's classification results 
the $B_{i} \otimes \Oh_{\infty}$ are isomorphic iff the $B_{i}$ are. Kirchberg has also 
shown that $\Oh_{2}$ absorbs any simple nuclear $C^{*}$-algebra tensorially, so in 
particular we have $B \otimes \Oh_{2} \cong \Oh_{2}$.  We will see later that $\Oh_{2}$, 
$\Oh_{\infty}$ and $B \otimes \Oh_{\infty}$ (with $B$ UHF of infinite type) are  the only 
purely infinite strongly self-absorbing $C^{*}$-algebras which satisfy the Universal 
Coefficients Theorem.  

\noindent
(iv) Let $p$, $q$ and $n$ be natural numbers with $p$ and $q$ dividing $n$.  
$C^*$-algebras of the form
\begin{displaymath}
I[p,n,q] = \{f \in M_n(\Ch([0,1])) \, | \, f(0)=1_{n/p} \otimes a, 
f(1)=b \otimes 1_{n/q},
a \in M_p, b \in M_q\}
\end{displaymath}
are commonly referred to as dimension drop intervals.  If $n=pq$ and $\mathrm{gcd}(p,q) 
= 1$, then the dimension drop
interval is said to be prime.  

In \cite{JS1}, Jiang and Su construct a $C^{*}$-algebra $\Zh$, which is the unique simple unital inductive
limit of dimension drop intervals having $\mathrm{K}_0 = \Z$, $\mathrm{K}_1 
= 0$ and a unique normalised
trace.  It is a limit of prime dimension drop intervals where the matrix 
dimensions tend to infinity,
and there is a unital embedding of any prime dimension drop interval into 
$\Zh$. Jiang and Su show that $\Zh$ is strongly self-absorbing; that $\Zh$ is $\mathrm{K}_{1}$-injective is established in \cite{J}.
It was shown in \cite{JS1} that a simple unital $C^{*}$-algebra absorbs $\Zh$ 
tensorially if it is AF or purely infinite. Therefore, tensoring our previous examples with $\Zh$ will not yield any new examples.
In fact, the entire list of examples provided here is closed under taking tensor products. 
\end{nexamples}
\en

\section{An intertwining argument}\label{inter}

Below we recall a result of Kirchberg and R{\o}rdam, based on Elliott's proof that $\mathcal{O}_2
\otimes \mathcal{O}_2 \cong \mathcal{O}_2$, 
which provides a characterisation of when a $C^{*}$-algebra $A$ is $\Dh$-stable 
($\Dh$ being strongly self-absorbing).  The statement involves the multiplier 
algebra $\Mh(A)$. However, to prove permanence properties of $\Dh$-stability, 
a slightly modified version of Theorem \ref{rordams-intertwining} (which is 
proved in a similar way but avoids use of the multiplier algebra) will be useful.

\bn
\label{corona}
\begin{nnotation}
For a $C^{*}$-algebra $A$ we denote by $\prod_{\N} A$ the $C^{*}$-algebra of bounded sequences 
over $\N$ with values in $A$;  the ideal of sequences converging to zero is denoted by 
$\bigoplus_{\N} A$. We shall write $\Qh (A)$ for the quotient $\prod_{\N} A / \bigoplus_{\N} A$. 
There is a canonical embedding $\iota: A \hookrightarrow \prod_{\N} A$, given by mapping $A$ 
to the subalgebra of constant sequences; $\iota$ clearly passes to an embedding of $A$ into 
$\Qh(A)$. For convenience, we will often omit the $\iota$ and simply identify $A$ with its 
image in $\prod_{\N}A$ or $\Qh(A)$, respectively.  If $B \subset A$ is a subalgebra, then 
$\Qh(A) \cap B'$ denotes the relative commutant of $B$ in $\Qh(A)$.
\end{nnotation}
\en

\bn
\label{rordams-intertwining}
\begin{ntheorem} 
{\rm (cf.\ \cite{K} and \cite{R1}, Theorem 7.2.2)} Let $A$ and $\Dh$ be separable $C^{*}$-algebras and 
suppose that $\Dh$ is unital and strongly self-absorbing. Then there is an isomorphism $\varphi: 
A \to A\otimes \Dh$ iff there is a unital $*$-homomorphism
\[
\varrho: \Dh \to \Qh (\Mh(A)) \cap A' \, .
\]
Moreover, in this case the maps $\varphi$ and $\id_{A} \otimes \be_{\Dh}$ are a.u.\ equivalent.  
\end{ntheorem}
\en

\bn
\label{nonunital-intertwining}
For some purposes (cf.\ the two subsequent sections) another version of the above will be useful. 
In the following result we have to ask $\Dh$ to be $\mathrm{K}_{1}$-injective (see \ref{unitary-homotopies}), 
which, for the known examples of strongly self-absorbing $C^{*}$-algebras, is no restriction (cf.\ \ref{absorbing-examples}).

\begin{ntheorem}
Let $A$ and $\Dh$ be separable $C^{*}$-algebras and suppose that $\Dh$ is unital, strongly 
self-absorbing and $\mathrm{K}_{1}$-injective, i.e., the canonical homomorphism  $\Uh(\Dh)/\Uh_{0}(\Dh) \to \mathrm{K}_{1}(\Dh)$ is injective.
Then there is an isomorphism $\varphi: A \to A\otimes \Dh$ iff there is a $*$-homomorphism
\[
\sigma: A \otimes \Dh \to \Qh (A)
\]
satisfying 
\[
\sigma(a \otimes \be_{\Dh}) = a \; \forall \, a \in A \, ,
\]
and in this case  $\varphi \au \id_{A} \otimes \be_{\Dh}$.  
\end{ntheorem}
\en

\bn
\label{unit-interval}
Before proving the theorem, we need some intermediate results.

\begin{nlemma}
Let $A$ be a separable $C^{*}$-algebra, sitting as an ideal in a separable unital $C^{*}$-algebra $B$. Then, there is a sequence of unital $*$-homomorphisms 
\[
\beta_{n}: \Ch([0,1]) \to A^{+} \subset B
\]
such that  the following hold:
\begin{itemize}
\item[(i)] $\beta_{n}(\Ch_{0}([0,1))) \subset A$ (regarding 
$\Ch_{0}([0,1))$ as subalgebra of $\Ch([0,1])$ in the canonical way)
\item[(ii)] $\|\beta_{n}(h)a - h(0) \cdot a\| \stackrel{n \to \infty}{\longrightarrow} 0 \; \forall \, h \in \Ch([0,1]), \, a \in A$
\item[(iii)] $\|[\beta_{n}(h) ,b ]\| \stackrel{n \to \infty}{\longrightarrow} 0 \; \forall \, h \in \Ch([0,1]), \, b \in B$.
\end{itemize}
\end{nlemma}

\begin{nproof}
First, choose an approximate unit $(e_{n})_{n \in \N}$ for $A$ which is quasicentral for $B$; we may assume the $e_{n}$ 
to be positive and normalized, moreover, we may assume that $e_{n} e_{n+1} = e_{n} \; \forall \, n \in \N$ (cf.\ \cite{Ped}, 3.12.16). Define continuous functions $g_{n}$ and $h_{n} \in \Ch([0,1])$ for $n \in \N$ by
\[
g_{n}(t) := \left\{ 
\begin{array}{ll}
1 & \mbox{for }  t =0 \\
0 & \mbox{for } t \in [\frac{1}{n}, 1]\\
\mathrm{linear} & \mathrm{elsewhere} 
\end{array}
\right.
\]
and
\[
h_{n}(t) := \left\{ 
\begin{array}{ll}
1 & \mbox{for }  t \in [0, \frac{1}{n}] \\
0 & \mbox{for } t =1\\
\mathrm{linear} & \mathrm{elsewhere} \, .
\end{array}
\right.
\]
One checks that for each $n$ there is a (unique) unital $*$-homomorphism 
\[
\beta_{n}: \Ch([0,1]) \to A^{+} 
\]
satisfying
\[
\beta_{n}(g_{n})= e_{n} \mbox{ and } \beta_{n}(h_{n})= e_{n+1} \, .
\]
Now if $h \in \Ch([0,1])$, then we have $\|h g_{n}  - h(0) \cdot g_{n}\| \to 0$, and so 
\[
\|\beta_{n}(h) e_{n} - h(0) \cdot e_{n}\| \to 0 \, ,
\] 
whence
\[
\|\beta_{n}(h)a - h(0) \cdot a\| \stackrel{n \to \infty}{\longrightarrow} 0 \; \forall \, a \in A\, .
\]
Moreover, there is a sequence of polynomials $(p_{k})_{k \in \N}$ in one variable, such that 
\[
\lim_{k \to \infty} \|p_{k}(\be - \id_{[0,1]}) -h\| =0 \, .
\]
But then we also have
\[
\lim_{k \to \infty} \lim_{n \to \infty} \|p_{k}(h_{n}) -h\| = 0 \, .
\]
Now since, for each $n$ and $k$, 
\[
\|p_{k}(e_{n+1}) - \beta_{n}(h)\| = \|\beta_{n}(p_{k}(h_{n}) - h)\| \le  \|p_{k}(h_{n}) -h\| \, ,
\]
we obtain
\[
\limsup_{k \to \infty} \lim_{n \to \infty} \|p_{k}(e_{n+1}) - \beta_{n}(h)\| = 0 \, .
\]
Note that 
\[
\lim_{n \to \infty} \|[p_{k}(e_{n+1}),b]\| = 0 \; \forall \, b \in B \, ,
\]
since the $e_{n}$ are quasicentral. As a consequence we see that, for all $b \in B$,
\begin{eqnarray*}
\lim_{n \to \infty} \|[\beta_{n}(h) , b]\| & \le & \limsup_{k \to \infty} \lim_{n \to \infty} (\|[p_{k}(e_{n+1}), b]\| +2 \cdot \|p_{k}(e_{n+1}) - \beta_{n}(h)\| ) \\
& \le & \limsup_{k \to \infty} \lim_{n \to \infty} \|[p_{k}(e_{n+1}), b]\|\\
& & +2 \cdot \limsup_{k \to \infty} \lim_{n \to \infty}\|p_{k}(e_{n+1}) - \beta_{n}(h)\| \\
& = & 0 \, .
\end{eqnarray*}
\end{nproof}
\en

\bn
\label{unitarization}
\begin{nlemma}
Let $A$ and $\Dh$ be separable $C^{*}$-algebras and suppose that $\Dh$ is unital, strongly 
self-absorbing and $\mathrm{K}_{1}$-injective. 
Then, there is a sequence $(s_{n})_{\N}$ of contractions in $A \otimes \Dh \otimes \Dh$ satisfying the following for all $a \in A, \, d \in \Dh$:
\begin{itemize}
\item[(i)] $\| [s_{n}, a \otimes \be_{\Dh} \otimes \be_{\Dh}] \| \to 0$ 
\item[(ii)] $\|s_{n}^{*} (a \otimes \be_{\Dh} \otimes d) s_{n} - a \otimes d \otimes \be_{\Dh} \| \to 0$
\item[(iii)] $\|s_{n}^{*}s_{n}(a \otimes \be_{\Dh} \otimes \be_{\Dh}) - a \otimes \be_{\Dh} \otimes \be_{\Dh} \|\to 0$
\item[(iv)] $s_{n} + \be - (s_{n}s_{n}^{*})^{\halb}$ is a unitary in $(A \otimes \Dh \otimes \Dh)^{+}$, where $\be$ denotes the unit of $(A \otimes \Dh \otimes \Dh)^{+}$.
\end{itemize}
\end{nlemma}

\begin{nproof}
For $i=0,2$ consider functions $h_{i} \in \Ch_{0}([0,1])$ defined by
\[
h_{i}(t) := \left\{ 
\begin{array}{ll}
1 & \mbox{for }  t \in [0, \frac{i}{3}] \\
0 & \mbox{for } t \in [\frac{i+1}{3}, 1] \\
\mathrm{linear} & \mathrm{elsewhere} \, .
\end{array}
\right.
\]
By Lemma \ref{unit-interval} there are $*$-homomorphisms $\beta_{n}: \Ch_{0}([0,1)) \to A$ 
such that 
\[
\|\beta_{n}(h_{i})a - a\| \stackrel{n \to \infty}{\longrightarrow} 0
\]
for $i=0,2$ and all $a \in A$. By hypothesis there are unitaries $v_{n} \in \Dh \otimes \Dh$ such that 
$\|v_{n}^{*}(\be_{\Dh} \otimes d) v_{n} - d \otimes \be_{\Dh}\| \to 0$ for all $d \in \Dh$. 
By Proposition \ref{unitary-homotopies}, each $v_{n}$ may be chosen to be homotopic (via unitaries) to $\be_{\Dh \otimes \Dh}$, so there are unitaries 
\[
u_{n} \in \Ch([0,1]) \otimes \Dh \otimes \Dh \cong \Ch([0,1], \Dh \otimes \Dh)
\]
such that 
\[
u_{n}(t) := \left\{ 
\begin{array}{ll}
v_{n} & \mbox{for }  t \in [0, \frac{1}{3}) \\
u_{n}(t) & \mbox{for } t \in [\frac{1}{3}, \frac{2}{3}) \\
\be_{\Dh \otimes \Dh} & t \in [\frac{2}{3},1] \, .
\end{array}
\right.
\]
We have elements 
\[
\tilde{u}_{n}:= (h_{2} \otimes \be_{\Dh \otimes \Dh}) u_{n} \in \Ch_{0}([0,1)) \otimes \Dh \otimes \Dh \, ,
\]
so we may define
\[
s_{n}:= \beta_{n} \otimes \id_{\Dh \otimes \Dh} (\tilde{u}_{n}) \in A \otimes \Dh \otimes \Dh \, .
\]
Note that 
\[
s_{n} (\beta_{n}(h_{0}) \otimes \be_{\Dh \otimes \Dh}) = (\beta_{n}(h_{0}) \otimes \be_{\Dh \otimes \Dh}) s_{n} = \beta_{n}(h_{0}) \otimes v_{n} \, .
\]
Together with
\[
\lim_{n \to \infty} \|\beta_{n}(h_{0})a - a\| = \lim_{n \to \infty} \|a \beta_{n}(h_{0})-a\| = 0 \; \forall \, a \in A
\]
this implies 
\[
\|(a \otimes x) s_{n} - a \otimes x v_{n} \| \stackrel{n \to \infty}{\longrightarrow} 0 \; \forall \, a \in A, \, x \in \Dh \otimes \Dh \, .
\]
But now it is straightforward to check that
\[
\|[s_{n}, a \otimes \be_{\Dh \otimes \Dh}]\| \to 0 \, , 
\]
that
\[
\|s_{n}^{*} (a \otimes \be_{\Dh} \otimes d) s_{n} - a \otimes d \otimes \be_{\Dh}\| \to 0
\]
and that 
\[
\|s_{n}^{*}s_{n} (a \otimes \be_{\Dh \otimes \Dh} ) - a \otimes \be_{\Dh \otimes \Dh}\| \to 0
\]
for all $a \in A, \, d \in \Dh$.

To check that $s_{n}+ \be - (s_{n}^{*}s_{n})^{\halb}$ is a unitary for each $n$, observe that
\[
(s_{n}^{*}s_{n})^{\halb}= \beta_{n}(h_{2}) \otimes \be_{\Dh \otimes \Dh} \, ,
\]
hence
\[
[(s_{n}^{*}s_{n})^{\halb}, s_{n}] = [(s_{n}^{*}s_{n})^{\halb}, s_{n}^{*}] = 0
\]
and
\[
s_{n}^{*} (\be - (s_{n}^{*}s_{n})^{\halb}) = (\be - (s_{n}^{*}s_{n})^{\halb}) s_{n} = \beta_{n}(h_{2} - h_{2}^{2}) \otimes \be_{\Dh\otimes \Dh} \, ,
\]
where for the last two identities we have used the definitions of $h_{2}$ and $u_{n}$. Now one computes 
\begin{eqnarray*}
\lefteqn{(s_{n}^{*} + \be - (s_{n}^{*}s_{n})^{\halb})(s_{n} + \be - (s_{n}^{*}s_{n})^{\halb})}\\
& = & \be + \beta_{n}(h_{2}^{2}+ 2 (h_{2}-h_{2}^{2})  -2 h_{2} + h_{2}^{2} ) \otimes \be_{\Dh \otimes \Dh} \\
& = & \be
\end{eqnarray*}
and, similarly,
\[
(s_{n} + \be - (s_{n}^{*}s_{n})^{\halb})(s_{n}^{*} + \be - (s_{n}^{*}s_{n})^{\halb}) = \be \, .
\]
\end{nproof}
\en

\bn
\label{corona-intertwining}
\begin{nprop}
Let $A$ and $B$ be separable $C^{*}$-algebras and $\iota: A \to B$ an embedding. Suppose 
there is a sequence of unitaries $v_{n} \in \Qh(B^{+})$ such that, for all $a \in A$ and $b\in B$,  $\|[v_{n}, \iota(a)]\| \to 0$ and $\dist(v_{n}^{*}bv_{n}, \Qh(\iota(A))) \to 0$.\\
Then, there is an isomorphism $\psi: A \to B$ which is a.u.\ equivalent to $\iota$.
\end{nprop}

\begin{nproof}
Let $\{a_{1}, \ldots , a_{k} \} \subset A$ and $\{b_{1}, \ldots, b_{l}\} \subset B$ be 
finite subsets of positive normalised elements and let $\varepsilon >0$.
In view of \cite{R1}, Proposition 2.3.5, it will suffice construct a unitary $u \in B^{+}$ 
such that $\|[u, \iota(a_{i})]\|< \varepsilon$ and $\dist(u^{*} b_{j} u, \iota(A)) < \varepsilon$ for $i=1, \ldots, k$ and $j=1, \ldots, l$.
By hypothesis, there are $N \in \N$ and $c_{1}, \ldots, c_{l} \in \Qh(\iota(A))$ 
satisfying $\|[v_{N}, \iota(a_{i})]\| < \varepsilon/2$, and $\|v_{N}^{*} b_{j} v_{N} - c_{j}\| < 
\varepsilon/2$ for $i=1, \ldots, k$ and $j=1, \ldots, l$. We may assume the $c_{j}$ to be positive and normalised as well.
Next, choose lifts $(u_{1}, u_{2}, \ldots) \in \prod B^{+}$ and $(\iota(a_{j,1}), \iota(a_{j,2}), \ldots ) 
\in \prod \iota(A) \subset \prod B$ for $v_{N}$ and $c_{j}$, $j=1, \ldots, l$, respectively. 
We may assume that the $u_{m}$ are unitaries (using \cite{R1}, Lemma 6.2.4) and (using functional calculus) that the $a_{j,m}$ are positive and normalised.
But since the $(u_{m})_{\N}$ and the $\iota(a_{j,m})_{m \in \N}$ are lifts for $v_{N}$ 
and $c_{j}$, there is $M \in \N$ such that $\|[u_{M}, \iota(a_{i})]\| < \varepsilon$ and 
$\| u_{M}^{*} b_{j} u_{M} - \iota(a_{j,M})\| < \varepsilon$ for $\iota = 1, \ldots, k$ and $j= 1, \ldots, l$.
\end{nproof}
\en

\begin{nproof}
(of Theorem \ref{nonunital-intertwining}) If $A \cong A \otimes \Dh$, for $n \in \N$ define $*$-homomorphisms 
\[
\sigma_{n}: A \otimes \Dh \otimes \Dh \to A \otimes \Dh
\]
by $\sigma_{n}:= \id_{A} \otimes \varphi_{n}$, where the $\varphi_{n}$ come from 
Proposition \ref{infinite-products}(iii).  The induced map $\sigma: (A \otimes \Dh) 
\otimes \Dh \to \Qh(A \otimes \Dh)$ obviously is a $*$-homomorphism satisfying 
\[
\sigma(a \otimes d \otimes \be_{\Dh})= a \otimes d \; \forall \, a \in A, \, d \in \Dh \, .
\]
  
Conversely, suppose there is a $*$-homomorphism $\sigma$ as in the Theorem. Let $\iota : A \to 
A \otimes \Dh$ be the canonical embedding given by $\id_{A} \otimes \be_{\Dh}$. Define 
$\beta: \Dh \to \Qh(A^{+} \otimes \Dh)$ by $\beta(d):= \be_{A^{+}} \otimes d$. Regarding 
$\Qh(A) \cong \Qh(\iota(A)) = \Qh(A \otimes \C \cdot \be_{\Dh})$ as subalgebra of $\Qh(A^{+} 
\otimes \Dh)$, we see that $\iota \circ \sigma$ and $\beta$ have commuting images and 
therefore induce a $*$-homomorphism
\[
\varrho: A \otimes \Dh \otimes \Dh \to \Qh(A^{+} \otimes \Dh)
\]
satisfying
\[
\varrho(a \otimes d_{0} \otimes d_{1}) = \iota \circ \sigma(a \otimes d_{0}) \beta(d_{1}) \, .
\]
Since $\iota \circ \sigma(A\otimes \Dh) \subset \Qh(A \otimes \Dh) \lhd \Qh(A^{+} \otimes \Dh)$, we see that in fact 
\[
\varrho(A \otimes \Dh \otimes \Dh) \subset \Qh(A \otimes \Dh) \subset \Qh(A^{+} \otimes \Dh) \, .
\]
By Lemma \ref{unitarization} there is a sequence $(s_{n})_{\N}$ of contractions in $A \otimes 
\Dh \otimes \Dh$ satisfying the following for all $a \in A$, $d \in \Dh$:
\begin{itemize}
\item[(i)] $\| [s_{n}, a \otimes \be_{\Dh} \otimes \be_{\Dh}] \| \to 0$ 
\item[(ii)] $\|s_{n}^{*} (a \otimes \be_{\Dh} \otimes d) s_{n} - a \otimes d \otimes \be_{\Dh} \| \to 0$
\item[(iii)] $\|s_{n}^{*}s_{n}(a \otimes \be_{\Dh} \otimes \be_{\Dh}) - a \otimes \be_{\Dh} \otimes \be_{\Dh} \|\to 0$
\item[(iv)] $s_{n} + \be - (s_{n}s_{n}^{*})^{\halb}$ is a unitary in $(A \otimes \Dh \otimes \Dh)^{+}$, 
where $\be$ again denotes the unit of $(A \otimes \Dh \otimes \Dh)^{+}$.
\end{itemize}
Set $v_{n}:= \varrho^{+}(s_{n}+ \be - (s_{n}^{*}s_{n})^{\halb})$, where 
\[
\varrho^{+}:(A \otimes \Dh \otimes \Dh)^{+} \to (\Qh(A \otimes \Dh))^{+} \subset \Qh((A \otimes \Dh)^{+})
\]
is the unitization of $\varrho$. Then 
\begin{eqnarray*}
\|[v_{n}, \iota(a)]\| & = & \|[v_{n}, \iota \circ \sigma(a \otimes \be_{\Dh}) \beta(\be_{\Dh})]\| \\
& = & \|[\varrho(s_{n} - (s_{n}^{*}s_{n})^{\halb}), \varrho(a \otimes \be_{\Dh \otimes \Dh})]\| \\
& = & \|\varrho([s_{n} - (s_{n}^{*}s_{n})^{\halb}, a \otimes \be_{\Dh \otimes \Dh}]) \| \\
& \stackrel{n \to \infty}{\longrightarrow} & 0 \, .
\end{eqnarray*}
Furthermore,
\begin{eqnarray*}
v_{n}^{*}(a \otimes d) v_{n} & = & v_{n}^{*} \iota \circ \sigma (a \otimes \be_{\Dh}) \beta(d) v_{n} \\
& = & v_{n}^{*} \varrho(a \otimes \be_{\Dh} \otimes d) v_{n} \\
& = & \varrho^{+}( v_{n}^{*} (a \otimes \be_{\Dh} \otimes d)v_{n}) \\
& \stackrel{n \to \infty}{\longrightarrow} & \varrho(a \otimes d \otimes \be_{\Dh}) \\
& = & \iota \circ \sigma (a \otimes d) \in \Qh(\iota(A)) \, ,
\end{eqnarray*}
so $\dist(v_{n}^{*}b v_{n}, \Qh(\iota(A))) \to 0 \; \forall \, b \in A \otimes \Dh$.

Finally, the $v_{n}$ are unitaries in $(\Qh(A \otimes \Dh))^{+} \subset \Qh((A \otimes \Dh)^{+})$, 
since $\varrho^{+}$ is a $*$-homomorphism, and the results follow from Proposition \ref{corona-intertwining}.
\end{nproof}

\bn
\label{embedding-products-2}
\begin{nremark}
The first part of the proof also shows that, if $A$ is $\Dh$-stable, there exists a 
sequence $(\sigma_{n}: A \otimes \Dh \to A)_{\N}$ of $*$-homomorphisms which  satisfies 
\[
\|\sigma_{n}(a \otimes \be_{\Dh}) - a\| \stackrel{n \to \infty}{\longrightarrow} 0 \; \forall \, a \in A \, .
\] 
\end{nremark}
\en

\section{Permanence properties of $\Dh$-stability}\label{perm}

In the sequel we conclude from Theorem \ref{nonunital-intertwining} that $\Dh$-stability passes to hereditary subalgebras 
(hence is Morita-invariant), quotients, inductive limits and to extensions. In the cases 
where $\Dh$ equals $\Oh_{2}$ or $\Oh_{\infty}$ these results were obtained independently 
(and with slightly different methods) by E.\ Kirchberg, cf.\ \cite{K2}, Section 8.  In \cite{HW}, I.\ Hirshberg and the second named author will show that $\Dh$-stability also passes to crossed products with Rokhlin actions of $\Z$, $\R$ or compact second countable groups.

Throughout this section, we shall assume $\Dh$ to be separable, unital, strongly self-absorbing 
and $\mathrm{K}_{1}$-injective (recall that the latter holds for all the examples of \ref{absorbing-examples}).

\bn
\label{hereditary-D-stable}
\begin{ncor}
Let $A$ be separable and $\Dh$-stable, and let $B \subset_{\her} A$ be a hereditary subalgebra. Then, $B$ is $\Dh$-stable.
\end{ncor}

\begin{nproof}
Let $\iota:B \to A$ be the injection map. Choose an  approximate unit $(h_{n})_{\N} \subset B$ for 
$B$; we assume the $h_{n}$ to be positive contractions. Let $h$ denote the image of the sequence
 $(h_{n})_{\N}$ in $\Qh(B) \subset \Qh(A)$. Furthermore, let $\beta: \Qh(A) \to \Qh(B)$ be 
the c.p.c.\ map given by $\beta(x):= h x h$. For $b \in B \subset A \subset \Qh(A)$ we have $hb=bh=b$; 
in particular, we obtain $\beta (b)=b \; \forall \, b \in B$.    Consider a $*$-homomorphism $\sigma: 
A \otimes \Dh \to \Qh(A)$ as in Theorem \ref{nonunital-intertwining} and define a c.p.c.\ map $\bar{\sigma}: 
B \otimes \Dh \to \Qh(B)$ by $\bar{\sigma}:= \beta \circ \sigma \circ (\iota \otimes \id_{\Dh})$. For $b \in B_{+}$ and $d \in \Dh_{+}$, we have 
\[
\bar{\sigma}(b \otimes \be_{\Dh})= \beta \circ \sigma (b \otimes \be_{\Dh}) = b
\]
and 
\begin{eqnarray*}
\bar{\sigma}(b \otimes d) & = & \beta ( \sigma(b^{\frac{1}{4}} \otimes \be_{\Dh})\sigma(b^{\halb} \otimes d)\sigma(b^{\frac{1}{4}} \otimes \be_{\Dh}) ) \\
& = & h(b^{\frac{1}{4}} \otimes \be_{\Dh})\sigma(b^{\halb} \otimes d)(b^{\frac{1}{4}} \otimes \be_{\Dh})h \\
& = & (b^{\frac{1}{4}} \otimes \be_{\Dh})\sigma(b^{\halb} \otimes d)(b^{\frac{1}{4}} \otimes \be_{\Dh})\\ 
& = & \sigma(b \otimes d) \, .
\end{eqnarray*}
The last equation not only shows that $\sigma$ maps $B \otimes \Dh$ to $\Qh(B) \subset \Qh(A)$, 
but also that $\bar{\sigma}$ is multiplicative, since $\sigma $ is. Therefore, $\bar{\sigma}$ 
satisfies the conditions of Theorem \ref{nonunital-intertwining}, hence $B$ is $\Dh$-stable.
\end{nproof}
\en

\bn
\label{morita-invariant}
\begin{ncor}
If $A$ is a  separable $C^{*}$-algebra and $r \in \N$, then $A$ is $\Dh$-stable iff 
$A \otimes M_{r}$ is $\Dh$-stable iff $A \otimes \Kh$ is $\Dh$-stable.   
\end{ncor}

\begin{nproof}
If $A$ is $\Dh$-stable, then clearly $A \otimes M_{r}$ and $A \otimes \Kh$ are. Conversely, 
if $A \otimes \Kh$ is $\Dh$-stable, then so is $A$ by Corollary \ref{hereditary-D-stable}.
\end{nproof}
\en

\bn
\label{quotients-D-stable}
\begin{ncor}
If $A$ is separable and $\Dh$-stable and $J \lhd A$ is a closed two-sided ideal, then $J$ and  $A/J$ are $\Dh$-stable. 
\end{ncor}

\begin{nproof}
$J \subset A$ is a hereditary subalgebra, so it is $\Dh$-stable by Corollary \ref{hereditary-D-stable}. For the second statement, note that the quotient map $q: A \to A/J$ induces a $*$-homomorphism 
$\bar{q}: \Qh(A) \to \Qh(A/J)$. Now let $\sigma: A \otimes \Dh \to 
\Qh(A)$ be a $*$-homomorphism as in Theorem \ref{nonunital-intertwining}. The 
composition $\bar{q} \circ \sigma$ passes to a $*$-homomorphism $\bar{\sigma}: (A/J) \otimes 
\Dh \to \Qh(A/J)$, because $\bar{q} \circ \sigma$ maps $J \otimes \be_{\Dh}$ to 
$0 \in \Qh(A/J)$, whence $\bar{q} \circ \sigma(J \otimes \Dh) = 0$. The map 
$\bar{\sigma}$ satisfies $\bar{\sigma} \circ (q \otimes \id_{\Dh}) = \bar{q} \circ \sigma$, hence 
\[
\bar{\sigma}(q(a) \otimes \be_{\Dh}) = \bar{q}(\sigma(a \otimes \be_{\Dh})) = \bar{q}(a) = q(a) \in \Qh(A/J)
\]  
for all $a \in A$. By surjectivity of $q$, the result now follows from Theorem \ref{nonunital-intertwining}. 
\end{nproof}
\en

\bn
\label{limits-D-stable-2}
\begin{ncor}
Let $A = \lim_{\to} A_{i}$ be an inductive limit of separable $\Dh$-stable $C^{*}$-algebras $A_{i}$, $i \in \N$. Then, $A$ is $\Dh$-stable. 
\end{ncor}

\begin{nproof}
Replacing the $A_{i}$ by their images in $A$ if necessary, by Corollary \ref{quotients-D-stable} 
we may assume the $A_{i}$ to form an increasing sequence of $\Dh$-stable $C^{*}$-algebras. 
From Remark \ref{embedding-products-2} for each $i \in \N$ we obtain $*$-homomorphisms
\[
\sigma_{i,n}: A_{i} \otimes \Dh \to A_{i} \subset A
\]
satisfying
\[
\sigma_{i,n}(a \otimes \be_{\Dh}) \stackrel{n \to \infty}{\longrightarrow} a \; \forall \, a \in A_{i} \,.
\]
Using separability of the $A_{i}$ we can find a sequence $(n_{i})_{i \in \N} \subset \N$ such that, for all $j \in \N$,
 \[
 \sigma_{i,n_{i}}(a \otimes \be_{\Dh}) \stackrel{i \to \infty}{\longrightarrow} a \; \forall \, a \in A_{j} \, . 
 \]
Note that the last statement makes sense even though $\sigma_{i,n_{i}}$ is only defined 
on $A_{j} \otimes \Dh$ for $j \le i$. Next we define a map $\tilde{\sigma}: \bigcup_{i} 
A_{i} \otimes \Dh \to \prod_{i} A_{i} \subset \prod_{i} A$ by
\[
\tilde{\sigma}_{i}(x) := \left\{
\begin{array}{ll}
\sigma_{i,n_{i}}(x) & \mbox{if } x \in A_{i} \otimes \Dh \\
0 & \mbox{else}.
\end{array}
\right.
\]
It is straightforward to see that the $\tilde{\sigma}$ induce a map 
\[
\bar{\sigma}: \bigcup_{i} A_{i} \otimes \Dh \to \Qh(A) 
\]
which is multiplicative, $*$-preserving and satisfies
\[
\bar{\sigma}(a \otimes \be_{\Dh})=a \; \forall \, a \in \bigcup_{i} A_{i} \, .
\]
Since $\bar{\sigma}$ is a $*$-homomorphism on  $A_{i} \otimes \Dh$, it is normdecreasing on $A_{i} \otimes \Dh$ for each $i \in \N$, hence on all of $\bigcup_{i} A_{i} \otimes \Dh$. Regarding $\bigcup_{i} A_{i} \otimes \Dh$ as a (dense) subalgebra of $A \otimes \Dh$, we see that $\bar{\sigma}$ extends to a $*$-homomorphism
\[
\sigma: A \otimes \Dh \to \Qh(A) \, ,
\]
still satisfying $\sigma(a \otimes \be_{\Dh})= a \; \forall \, a \in A$, whence $A$ is $\Dh$-stable by Theorem \ref{nonunital-intertwining}.  
\end{nproof}
\en

\section{Extensions}\label{ext}

We have already seen that $\Dh$-stability passes to quotients and ideals; in this section we show that it is also stable under taking extensions.

\bn
\label{multiplicative-domain}
In the proof of Theorem \ref{extensions-D-stable} below we shall have use for  a straightforward and well-known 
consequence of Stinespring's theorem (cf.\ \cite{KW},
Lemma 3.5):

\begin{nlemma}
Let $A$, $B$ be $C^*$-algebras and $\varphi: B \to A$ a c.p.c.\ map. Then, 
for any $x,y \in B_{+}$, we have $\| \varphi(xy)
- \varphi(x) \varphi(y)\| \le \|\varphi(x^{2}) - 
\varphi(x)^{2}\|^{\halb}  \|y\|$.
In particular, if $x$ is in the multiplicative domain of $\varphi$, i.e., $\|\varphi(x^{2}) - 
\varphi(x)^{2}\|=0$, then $\varphi(xy)=\varphi(x) \varphi(y)$ for all $y \in B_{+}$, hence for all $y \in B$. 
\end{nlemma}
\en

\bn
\label{commutant-onto}
The next result is only a minor variation of \cite{K2}, Lemma 2.6(iii).

\begin{nlemma}
Let $0 \to J \to E \to A \stackrel{q}{\to} 0$ be a short exact sequence of separable 
$C^{*}$-algebras. Then the induced map $\bar{q}: \Qh(E) \to \Qh(A)$ maps $\Qh(E) \cap E'$ onto $\Qh(A) \cap A'$.
\end{nlemma}

\begin{nproof}
Since $\Qh(A) \cap A'$ is a $C^{*}$-algebra, it suffices to show that any positive 
contractive $a \in \Qh(A) \cap A'$ lifts to some $e \in \Qh(E) \cap E'$. So let $a$ 
be represented by a sequence $(a_{n})_{n \in \N}$ of positive contractions in $A$ satisfying 
\[
\|[a_{n},x]\| \stackrel{n \to \infty}{\longrightarrow} 0 \; \forall \, x \in A \, .
\] 
Each $a_{n}$ lifts to a positive contraction $e_{n} \in E$. Now choose a quasicentral 
approximate unit $(d_{n})_{n \in \N}$ for $J$ and a sequence $(k_{n})_{n \in \N} \subset \N$ such that
\[
\|[y,(\be_{E^{+}}-d_{k_{n}})e_{n}]\| \stackrel{n \to \infty}{\longrightarrow} 0 \; \forall \, y \in E \, ;
\]
this is possible since $(d_{k})_{k \in \N}$ is quasicentral with respect to $E$ and 
$q([y,e_{n}]) \stackrel{n \to \infty}{\longrightarrow} 0$, whence $\dist([y,e_{n}],J) \stackrel{n \to \infty}{\longrightarrow} 0$. 
Let $e \in \Qh(E)$ be the element represented by $((\be_{E^{+}}-d_{k_{n}})e_{n})_{n \in \N}$, then $e \in \Qh(E)\cap E'$ and 
\[
\bar{q}(e)= [(q((\be_{E^{+}}-d_{k_{n}})e_{n}))_{n \in \N}] = [(q(e_{n}))_{n \in \N}] = [(a_{n})_{n \in \N}] = a \, .
\] 
\end{nproof}
\en

\bn
\label{extensions-D-stable}
\begin{ntheorem}
Let $0 \to J \to E \stackrel{q}{\to} A \to 0$ be a short exact sequence of separable 
$C^{*}$-algebras. With $\Dh$ as in the preceding section, if $J$ and $A$ are $\Dh$-stable,  so is $E$.
\end{ntheorem}

\begin{nproof}
The proof is quite technical, so we briefly sketch its idea. We want to construct a $*$-homomorphism 
\[
\gamma: E \otimes \Dh \to \Qh(E)
\]
such that $\gamma|_{E \otimes \be_{\Dh}} = \id_{E}$. $\gamma$ will be a refined convex combination of c.p.c.\ maps 
\[
\bar{\varrho}, \, \mu : E \otimes \Dh \to \Qh(E) \, , 
\]
where $\bar{\varrho}$ and $\mu$ are constructed from the $*$-homomorphisms $A \otimes \Dh \to \Qh(A)$ and $J \otimes \Dh \to \Qh(J)$ implementing $\Dh$-stability of $A$ and $J$, respectively. At the same time, a quasicentral approximate unit of $J$ will yield a unital $*$-homomorphism 
\[
\beta: \Ch([0,1]) \to \Qh(E^{+}) \, .
\]
The map $\gamma$ is a combination of $\bar{\varrho}$ and $\mu$ `along' the image of the unit interval under $\beta$. On the left hand side of the interval, $\gamma$ coincides with $\bar{\varrho}$, on the right hand side with $\mu$. The problem is, that $\bar{\varrho}$ and $\mu$ yield two distinct copies of $\Dh$. However, using the assumptions on $\Dh$, we can construct a continuous path of unitaries which connects the unit of $\Dh \otimes \Dh$ with a unitary implementing the half-flip on $\Dh \otimes \Dh$. This path is then used to intertwine the two copies of $\Dh$ along the middle part of the interval. If all these maps are chosen carefully enough, we then obtain the desired $*$-homomorphism $\gamma$. All this will now be made precise.

Let $(\varrho_{n}: A \otimes \Dh \to A)_{n \in \N}$  be a sequence of $*$-homomorphisms  satisfying 
\begin{equation}
\label{almost-central}
\|\varrho_{n}(a \otimes \be_{\Dh}) - a\| \stackrel{n \to \infty}{\longrightarrow} 0 \; \forall \, a \in A \,  ;  
\end{equation}
such a sequence exists by Remark \ref{embedding-products-2}. The $\varrho_{n}$ induce a $*$-homomorphism 
\[
\varrho: A \otimes \Dh \to \Qh(A) \, .
\]
Let $(h_{n})_{\N}$ be an approximate unit for $A$. Define c.p.c.\ maps $\varrho'_{n}: 
\Dh \to A$ by $\varrho'_{n}(d):= \varrho_{n}(h_{n} \otimes d)$. The induced c.p.c.\ map 
$\varrho': \Dh \to \Qh(A)$ in fact maps $\Dh$ to $\Qh(A) \cap A'$. Moreover, one checks that 
\begin{equation}
\label{varrhoprime}
\varrho(a \otimes d)= a \varrho'(d) = \varrho'(d)a \; \forall \, a \in A, \, d \in \Dh \, .
\end{equation}
Let $\bar{q}: \Qh(E) \to \Qh(A)$ denote the obvious map induced by $q$. Since $\Dh$ is nuclear, Lemma \ref{commutant-onto} and the Choi--Effros lifting theorem imply that $\varrho'$ has  a c.p.c.\ lift 
\[
\tilde{\varrho}: \Dh \to \Qh(E) \cap E' \, ,
\]
which in turn lifts to a sequence of c.p.c.\ maps $\tilde{\varrho}_{n}: \Dh \to E$. Since the image of $\tilde{\varrho}$ commutes with $E$, there is a c.p.c.\ map 
\[
\bar{\varrho}: E \otimes \Dh \to \Qh(E) \, ,
\]
given by 
\[
\bar{\varrho}(e \otimes d) = e \tilde{\varrho}(d) = \tilde{\varrho}(d) e \; \forall \, e \in E, \, d \in \Dh \, . 
\]
Since
\begin{eqnarray*}
\bar{q} \circ \bar{\varrho}(e \otimes d) & = & q(e) \varrho'(d) \\
& = & \varrho(q(e) \otimes d) \\
& = & \varrho \circ (q \otimes \id_{\Dh})(e \otimes d) \, ,
\end{eqnarray*}
we have 
\[
\varrho \circ (q \otimes \id_{\Dh})= \bar{q}  \circ \bar{\varrho} \, .
\]
Together with (\ref{almost-central}) and (\ref{varrhoprime}) this shows that 
\begin{equation}
\label{q-varrho-unit}
\|q(e \tilde{\varrho}_{n}(\be_{\Dh}) - e)\| \to 0 \; \forall \, e \in E_{+}
\end{equation}
and that $\bar{q} \circ \bar{\varrho}$ is a $*$-homomorphism, whence 
\begin{equation}
\label{q-multiplicative}
\|q(e \tilde{\varrho}_{n}(d^{2})- e \tilde{\varrho}_{n}(d)^{2})\| \to 0  \; \forall \, e \in E_{+}, \, d \in \Dh_{+} \, .
\end{equation}

Before proceeding, we define continuous functions on the unit interval as follows:
\[
g'_{0}(t):= \left\{ 
\begin{array}{ll}
1 & \mbox{for }t =0\\
0 & \mbox{for }t \ge \frac{1}{8}\\
\mbox{linear} & \mbox{else} \, ,
\end{array}
\right. 
\]
\[
g_{0}(t):= \left\{ 
\begin{array}{ll}
1 & \mbox{for }  t \le \frac{1}{8}\\
0 & \mbox{for }t \ge \frac{2}{8}\\
\mbox{linear} & \mbox{else} \, ,
\end{array}
\right. 
\]
and 
\[
g'_{1}(t):= g'_{0}(1-t), \, g_{1}(t):= g_{0}(1-t) , \, g_{\halb}:= 1 - g_{0} - g_{1} \, . 
\]
Applying Lemma \ref{unit-interval} (with $J$ in place of $A$ and $E^{+}$ in place of $B$) we may use a diagonal sequence argument to obtain unital 
$*$-homomorphisms 
\[
\beta_{n}: \Ch([0,1]) \to J^{+} \subset E^{+} \, ,
\]
$n \in \N$, with the following properties:
\begin{itemize}
\item[a)] $\|\beta_{n}(g_{0}')c - c\| \to 0 \; \forall \, c \in J$
\item[b)] $\|\beta_{n}(\be_{[0,1]} - g'_{0})(e \tilde{\varrho}_{n}( d^{2}) - e \tilde{\varrho}_{n}(d)^{2})\| \to 0 \; \forall \, e \in E_{+}, \, d \in \Dh_{+}$ (using (\ref{q-multiplicative}))
\item[c)] $\|\beta_{n}(\be_{[0,1]} - g'_{0})(e \tilde{\varrho}_{n}(\be_{\Dh}) - e)\| \to 0 \; \forall \, e \in E_{+}$ (using (\ref{q-varrho-unit}))
\item[d)] $\|[\beta_{n}(f), e]\| \to 0$ and $\|[\beta_{n}(f), e\tilde{\varrho}_{n}(d)]\| \to 0 \; \forall \, f \in \Ch([0,1]), \, e \in E, \, d \in \Dh $ (regarding $J^{+}$ as a subalgebra of $E^{+}$)
\item[e)] $\beta_{n}(f) \subset J \; \forall \, n \in \N, \, f \in \Ch_{0}([0,1))$, in particular $\be_{E^{+}}- \beta_{n}(g_{1}')  \in J \; \forall \, n \in \N$.
\end{itemize}
The $\beta_{n}$ induce a $*$-homomorphism $\beta: \Ch([0,1]) \to \Qh(J^{+}) \subset \Qh(E^{+})$ satisfying 
\[
\beta(\Ch([0,1])) \subset \bar{\varrho}(E \otimes \Dh)' \, , \, \beta(g_{0}')c = c   
\]
and
\begin{equation}
\label{-1}
\beta(\be_{[0,1]}- g'_{0})(\bar{\varrho}(e \otimes \be_{\Dh})- e) = \beta(\be_{[0,1]}- g'_{0})(\bar{\varrho}(e^{2} \otimes d^{2})- \bar{\varrho}(e \otimes d)^{2})= 0 
\end{equation} 
for all $c \in J_{+}$, $e \in E_{+}$ and $d \in \Dh_{+}$.

From Remark \ref{embedding-products-2}, we obtain $*$-homomorphisms $\zeta_{n}: J \otimes \Dh \to J$ satisfying
\begin{equation}
\label{0}
\|\zeta_{n}(c \otimes \be_{\Dh}) - c\| \to 0  \; \forall \, c \in J \, .
\end{equation}
With a little extra effort, using  (\ref{0}), e) and separability of $E$, we may even assume that
\begin{eqnarray}
\label{1}
\lefteqn{\| \zeta_{n}(((\be_{E^{+}} - \beta_{n}(g_{1}'))^{\halb} e \tilde{\varrho}_{n}(d) (\be_{E^{+}} - \beta_{n}(g_{1}'))^{\halb} ) \otimes \be_{\Dh} )} \nonumber\\
&& - (\be_{E^{+}} - \beta_{n}(g_{1}'))^{\halb} e \tilde{\varrho}_{n}(d) (\be_{E^{+}} - \beta_{n}(g_{1}'))^{\halb}  \| \to 0 \, ,
\end{eqnarray}
that
\begin{equation}
\label{1.5}
\| \zeta_{n} (\beta_{n}(g_{0}) e \otimes \be_{\Dh}) - \beta_{n}(g_{0})e \| \to 0 
\end{equation}
and that
\begin{equation}
\label{2}
\| \zeta_{n}(\beta_{n}(f) \otimes \be_{\Dh}) - \beta_{n}(f) \| \to 0 
\end{equation}
for all $e \in E^{+}$, $d \in \Dh$ and $f \in \Ch_{0}([0,1))$. 

Define $\mu_{n}: E^{+} \otimes \Dh \to J$ by
\[
\mu_{n}(x):= \zeta_{n} (((\be_{E^{+}} - \beta_{n}(g_{1}'))^{\halb} \otimes \be_{\Dh} )x ((\be_{E^{+}} - \beta_{n}(g_{1}'))^{\halb} \otimes \be_{\Dh} )) \, ; 
\]
the $\mu_{n}$ are well-defined by e); they are c.p.c.\ and one checks that the induced map $\mu: E^{+} \otimes \Dh \to \Qh(J)$ satisfies the following:
\begin{itemize}
\item[g)] $\mu|_{J \otimes \be_{\Dh}} = \id_{J}$ (by a), (\ref{0}) and the definition of the $g_{i}'$)
\item[h)] $\mu|_{J \otimes \Dh}$ is a $*$-homomorphism (by a), using that the $\zeta_{n}$ are $*$-homomorphisms)
\item[i)] $\range (\mu) \subset (\range (\beta))'$ (using d), (\ref{2}) and multiplicativity of the $\zeta_{n}$)
\item[j)] $(\be_{E^{+}} - \beta(g_{1}))(\mu(e^{2} \otimes d^{2}) - \mu(e \otimes d)^{2}) = 0 \; \forall \, e \in E_{+}, \, d \in \Dh_{+}$ (by (\ref{2}) and multiplicativity of the $\zeta_{n}$)
\item[k)] $\mu(e \otimes \be_{\Dh}) = \beta(\be_{[0,1]} - g_{1}') \bar{\varrho}(e \otimes \be_{\Dh}) \; \forall \, e \in E$ (by d) and (\ref{1}))
\item[l)] $\beta(g'_{0})x = x \; \forall \, x \in \mu(J \otimes \Dh)$ (by a), g) and h))
\item[m)] $\beta(\be_{[0,1]}- g'_{0}) \perp \mu(J \otimes \Dh)$ (also by a), g) and h))
\item[n)] $\| [\mu(\be_{E^{+}} \otimes d_{0}), (\be_{E^{+}} - \beta(g_{1}'))^{\halb} 
\bar{\varrho}(e \otimes d_{1})(\be_{E^{+}} - \beta(g_{1}'))^{\halb}] \| = 0 \; \forall \, e \in E, \, d_{i} \in \Dh$ (by (\ref{1}))
\item[o)] $\mu(\be_{E^{+}} \otimes \be_{\Dh}) = \be_{E^{+}} - \beta(g_{1}')$ (by (\ref{2}))
\item[p)] $\beta(g_{0}) \mu(e \otimes \be_{\Dh}) = \beta(g_{0}) e \; \forall \, e \in E$ (by (\ref{1.5}) and (\ref{2})).
\end{itemize}
In particular we see from (\ref{-1}), i) and n) that the c.p.c.\ maps 
\[
\beta:\Ch([0,1]) \to \Qh(E^{+}) \, , 
\]
\[
\ad((\be_{E^{+}} - \beta(g_{1}'))^{\halb}) \circ \bar{\varrho}:E \otimes \Dh \to \Qh(E)
\]
and 
\[
\mu: \C \cdot \be_{E^{+}} \otimes \Dh \to \Qh(J) \subset \Qh(E)
\]
have commuting images in $\Qh(E^{+})$ and thus give rise to a c.p.c.\ map
\[ \textstyle
\lambda: \Ch_{0}((\frac{1}{8},\frac{7}{8})) \otimes E \otimes \Dh \otimes \Dh \to \Qh(E)
\] 
satisfying
\[
\lambda(f \otimes e \otimes d_{0} \otimes d_{1}) = \beta(f) \cdot \mu(\be_{E^{+}} \otimes d_{0}) 
\cdot (\be_{E^{+}} - \beta(g_{1}'))^{\halb} \bar{\varrho}(e \otimes d_{1})(\be_{E^{+}} - \beta(g_{1}'))^{\halb} \, .
\]
Since $f= (\be_{[0,1]}-g_{1}')f \; \forall \, f \in \Ch_{0}((\frac{1}{8},\frac{7}{8}))$, we in fact have
\begin{equation}
\label{3}
\lambda(f \otimes e \otimes d_{0} \otimes d_{1}) = \beta(f) \cdot \mu(\be_{E^{+}} \otimes d_{0}) \cdot  \bar{\varrho}(e \otimes d_{1}) 
\end{equation}
and, using i), j) and Lemma \ref{multiplicative-domain}, that $\lambda$ is a $*$-homomorphism. Note that 
\begin{equation}
\label{3.5}
\lambda(f \otimes e \otimes \be_{\Dh} \otimes \be_{\Dh}) = \beta(f) \cdot \bar{\varrho}(e \otimes \be_{\Dh})
\end{equation}
(by o)) and that the image of $\lambda$ in fact lies in $\Qh(J)$, since the image of $\mu$ does and $\Qh(J) \lhd \Qh(E)$ is an ideal. 

Choose unitaries $s_{l} \in \Dh \otimes \Dh$ such that $\|s_{l}^{*}(x \otimes \be_{\Dh})s_{l} - \be_{\Dh} 
\otimes x\| \stackrel{l \to \infty}{\longrightarrow} 0 \; \forall \, x \in \Dh$. By our assumption on 
$\Dh$ (in connection with Proposition \ref{unitary-homotopies}) we may assume the $s_{l}$ to be homotopic 
to $\be_{\Dh \otimes \Dh}$; therefore, there are unitaries $\tilde{s}_{l} \in \Ch([0,1], \Dh \otimes \Dh)$ such that 
\[
\tilde{s}_{l}|_{[0, \frac{1}{4}]} \equiv \be_{\Dh \otimes \Dh} \mbox{ and } \tilde{s}_{l}|_{[\frac{3}{4},1]} \equiv s_{l} \,.
\]
We regard the $\tilde{s}_{l}$ as elements of $\Ch([0,1]) \otimes E^{+} \otimes \Dh \otimes \Dh$. 

Let $(d_{l})_{\N}$ be an approximate unit for $E$ (the $d_{l}$ being positive contractions) and define 
$v_{l} \in \Ch_{0}((\frac{1}{8},\frac{7}{8})) \otimes E \otimes \Dh \otimes \Dh$ by
\[
v_{l} := (g_{\halb}^{\frac{1}{4}} \otimes d_{l} \otimes \be_{\Dh} \otimes \be_{\Dh}) \cdot \tilde{s}_{l} \, .
\]
One checks that 
\[
v_{l}^{*}v_{l} = v_{l}v_{l}^{*} = g_{\halb}^{\halb} \otimes d_{l}^{2} \otimes \be_{\Dh} \otimes \be_{\Dh}
\]
and that, for any $f \in \Ch_{0}((\frac{1}{8},\frac{7}{8})), \, e \in E$ and $ x,y \in \Dh$, the sequences
\[
\begin{array}{l}
v_{l}^{*}(f \cdot g_{0} \otimes e \otimes x \otimes y) v_{l} -  g_{\halb}^{\halb} \cdot f \cdot g_{0} \otimes e \otimes x \otimes y \, ,\\
v_{l}^{*}(f \cdot g_{1} \otimes e \otimes x \otimes \be_{\Dh}) v_{l} -  g_{\halb}^{\halb} \cdot f \cdot g_{1} \otimes e \otimes \be_{\Dh} \otimes x \, ,\\
v_{l}^{*}(f  \otimes e \otimes \be_{\Dh} \otimes \be_{\Dh}) v_{l} -  g_{\halb}^{\halb} \cdot f  \otimes e \otimes \be_{\Dh} \otimes \be_{\Dh} \, \\
\end{array}
\]
all converge to zero as $l $ goes to infinity. 

For each $l \in \N$,  $\lambda(v_{l})$ is a contraction in $\Qh(J)$,  so there is a contractive 
lift $(v_{l,n})_{n\in \N} \in \prod_{\N} J$. For a suitable (in a sense to be made precise shortly) 
increasing sequence $(n_{l})_{l \in \N} \subset \N$ define 
\[
u:= [(0,\ldots,0,v_{0,n_{0}},v_{0,n_{0}+1}, \ldots, v_{0,n_{1}-1},v_{1,n_{1}}, \ldots, v_{1,n_{2}-1},v_{2,n_{2}}, \ldots)] \in \Qh(J) \, .
\]
It is not hard to see that $(n_{l})_{l\in \N}$ can be chosen such that the following hold for all $f \in \Ch_{0}((\frac{1}{8},\frac{7}{8}))$, $e \in E$ and  $x,y \in \Dh$: 
\begin{itemize}
\item[q)] $u^{*}u \lambda(f \otimes e \otimes x \otimes y) = u u^{*} \lambda(f \otimes e \otimes x \otimes y)\\
= \lambda(f \otimes e \otimes x \otimes y) u^{*} u = \lambda(f \otimes e \otimes x \otimes y) u u^{*} \\
= \lambda(g_{\halb}^{\halb} \cdot f \otimes e \otimes x \otimes y)$
\item[r)] $u^{*} \lambda(g_{0} \cdot f \otimes e \otimes x \otimes y) u = \lambda(g_{0} \cdot g_{\halb}^{\halb} \cdot f \otimes e \otimes x \otimes y)$
\item[s)] $u^{*} \lambda(g_{1} \cdot f \otimes e \otimes x \otimes \be_{\Dh}) u = \lambda(g_{1} \cdot g_{\halb}^{\halb} \cdot f \otimes e \otimes \be_{\Dh} \otimes x)$
\item[t)] $u^{*} \lambda(f \otimes e \otimes \be_{\Dh} \otimes \be_{\Dh})u = \lambda(g_{\halb}^{\halb} \cdot f \otimes e \otimes \be_{\Dh} \otimes \be_{\Dh})$.
\end{itemize}
Moreover, we have $u \in \bar{\varrho}(E \otimes \be_{\Dh})' \cap (\range (\beta))' $.

We note some simple computations for later use. For $e \in E_{+}$ and $d \in \Dh_{+}$ we have
\begin{eqnarray}
\label{5}
\lefteqn{\beta(g_{0}) \cdot \mu(e \otimes d) u^{*}\lambda(g_{\halb}^{\halb} \otimes e \otimes d \otimes \be_{\Dh})u} \nonumber \\
&=& \mu(e \otimes d) \lambda(g_{0} \cdot g_{\halb} \otimes e \otimes d \otimes \be_{\Dh}) \nonumber \\
&=& \mu(e \otimes d)  \beta(g_{0} \cdot g_{\halb}) \mu(\be_{E^{+}} \otimes d) \bar{\varrho}(e \otimes \be_{\Dh}) \nonumber \\
&\stackrel{\mathrm{k)}}{=}& \beta(g_{\halb})\beta(g_{0}) \mu(e \otimes d) \mu(\be_{E^{+}} \otimes d) \mu(e \otimes \be_{\Dh}) \nonumber \\
&\stackrel{\mathrm{j)}, \ref{multiplicative-domain}}{=}& \beta(g_{\halb})\beta(g_{0}) \mu(e \otimes d^{2})  \mu(e \otimes \be_{\Dh}) \nonumber \\ 
&\stackrel{\mathrm{j)}, \ref{multiplicative-domain}}{=}& \beta(g_{\halb})\beta(g_{0}) \mu(\be_{E^{+}} \otimes d^{2})  \mu(e \otimes \be_{\Dh})^{2} \nonumber \\ 
&\stackrel{\mathrm{k)}}{=}& \beta(g_{\halb})\beta(g_{0}) \mu(\be_{E^{+}} \otimes d^{2}) \bar{\varrho}(e^{2} \otimes \be_{\Dh}) \nonumber \\
&=& \lambda(g_{0} \cdot g_{\halb} \otimes e^{2} \otimes d^{2} \otimes \be_{\Dh}) \nonumber \\
&\stackrel{\mathrm{r)}}{=}& \beta(g_{0})u^{*} \lambda(g_{\halb}^{\halb} \otimes e^{2} \otimes d^{2} \otimes \be_{\Dh})u \\
&=& \beta(g_{\halb})\beta(g_{0}) \mu(e^{2} \otimes d^{2}) \, , \nonumber 
\end{eqnarray}
and, similarly,
\begin{eqnarray} 
\label{6}
 u^{*}\lambda(g_{\halb}^{\halb} \otimes e \otimes d \otimes \be_{\Dh})u \cdot  \beta(g_{0}) 
\cdot \mu(e \otimes d) & = & \beta(g_{\halb})\beta(g_{0}) \mu(e^{2} \otimes d^{2}) \, .
\end{eqnarray}
Furthermore, we see that
\begin{eqnarray}
\label{7}
\lefteqn{\beta(g_{1}) \bar{\varrho}(e \otimes d) u^{*}\lambda(g_{\halb}^{\halb} \otimes e \otimes d \otimes \be_{\Dh})u} \nonumber \\
& = & \bar{\varrho}(e \otimes d) \lambda(g_{1} \cdot g_{\halb} \otimes e \otimes \be_{\Dh} \otimes d) \nonumber \\
& \stackrel{\mathrm{o)}}{=} & \bar{\varrho}(e \otimes d) \beta(g_{1} \cdot g_{\halb}) \beta(\be_{[0,1]}- g_{1}')\bar{\varrho}(e \otimes d) \nonumber \\
& \stackrel{(\ref{-1})}{=} & \beta(g_{1}) \beta(g_{\halb}) \bar{\varrho}(e^{2} \otimes d^{2})  \\
& = & \beta(g_{1}) \lambda(g_{\halb} \otimes e^{2} \otimes \be_{\Dh} \otimes d^{2}) \nonumber \\
& = & \beta(g_{1}) u^{*} \lambda(g_{\halb}^{\halb} \otimes e^{2} \otimes d^{2} \otimes \be_{\Dh})u \nonumber
\end{eqnarray}
and that
\begin{equation}
\label{7.5}
u^{*}\lambda(g_{\halb}^{\halb} \otimes e \otimes d \otimes \be_{\Dh})u  \, \beta(g_{1}) 
\bar{\varrho}(e \otimes d) =  \beta(g_{1}) u^{*} \lambda(g_{\halb}^{\halb} \otimes e^{2} \otimes d^{2} \otimes \be_{\Dh})u \, .
 \end{equation}

Now define a c.p.c.\ map 
\[
\gamma: E \otimes \Dh \to \Qh(E)
\]
by
\[
\gamma(e \otimes d):= \beta(g_{0}) \cdot \mu(e \otimes d) + u^{*}\lambda(g_{\halb}^{\halb} \otimes 
e \otimes d \otimes \be_{\Dh})u + \beta(g_{1}) \cdot \bar{\varrho}(e \otimes d) \, .
\]
We proceed to check that $\gamma$ satisfies the conditions of Theorem \ref{nonunital-intertwining}, 
i.e., it is a $*$-homomorphism sending $e \otimes \be_{\Dh}$ to $e$ for all $e \in E$.

First, we compute
\begin{eqnarray*}
\gamma(e \otimes \be_{\Dh}) & = & \beta(g_{0}) \cdot \mu(e \otimes \be_{\Dh}) \\
&&+ u^{*} \lambda(g_{\halb}^{\halb} \otimes e \otimes \be_{\Dh} \otimes \be_{\Dh})u \\
&&+ \beta(g_{1}) \cdot \bar{\varrho}(e \otimes \be_{\Dh}) \\
& \stackrel{\mathrm{p),t)}}{=} & \beta(g_{0}) \cdot e \\
&&+ \lambda(g_{\halb} \otimes e \otimes \be_{\Dh} \otimes \be_{\Dh})\\
&& + \beta(g_{1}) \cdot \bar{\varrho}(e \otimes \be_{\Dh})\\
& \stackrel{(\ref{3.5})}{=} & \beta(g_{0}) \cdot e + \beta(g_{\halb}) \cdot \bar{\varrho}(e \otimes \be_{\Dh}) + \beta(g_{1}) \cdot \bar{\varrho}(e \otimes \be_{\Dh}) \\
&\stackrel{(\ref{-1})}{=}& \beta(g_{0}) \cdot e + \beta(g_{\halb}) \cdot e + \beta(g_{1}) \cdot e\\
&=& e \, . 
\end{eqnarray*}

Finally, we check that $\gamma$ is multiplicative on $E \otimes \Dh$. By Lemma \ref{multiplicative-domain}, 
we only have to show that $\gamma(e \otimes d)^{2} = \gamma(e^{2} \otimes d^{2})$ for $e \in E_{+}$ and $d \in \Dh_{+}$:
\begin{eqnarray*}
\lefteqn{\gamma(e \otimes d)^{2}}\\
&=& \beta(g_{0})^{2} \cdot \mu(e \otimes d)^{2} + u^{*} \lambda(g_{\halb}^{\halb} \otimes e \otimes d \otimes 
\be_{\Dh}) u u^{*} \lambda(g_{\halb}^{\halb} \otimes e \otimes d \otimes \be_{\Dh}) u \\
& & + \beta(g_{1})^{2} \cdot \bar{\varrho}(e \otimes d)^{2} + \beta(g_{0}) \cdot \mu(e \otimes d) u^{*} 
\lambda(g_{\halb}^{\halb} \otimes e \otimes d \otimes \be_{\Dh}) u \\
& & +u^{*} \lambda(g_{\halb}^{\halb} \otimes e \otimes d \otimes \be_{\Dh}) u \beta(g_{0}) \mu(e \otimes d) \\
&&+ \beta(g_{1}) \bar{\varrho}(e \otimes d) u^{*} \lambda(g_{\halb}^{\halb} \otimes e \otimes d \otimes \be_{\Dh}) u \\
& & +u^{*} \lambda(g_{\halb}^{\halb} \otimes e \otimes d \otimes \be_{\Dh}) u \beta(g_{1}) \bar{\varrho}(e \otimes d)\\
& = & 
\beta(g_{0})^{2} \cdot \mu(e^{2} \otimes d^{2})  + \beta(g_{\halb}) \cdot u^{*} \lambda(g_{\halb}^{\halb} \otimes e^{2} \otimes d^{2} \otimes \be_{\Dh}) u \\
&& + \beta(g_{1})^{2} \cdot \bar{\varrho}(e^{2} \otimes d^{2}) + \beta(g_{0}) \cdot u^{*} \lambda(g_{\halb}^{\halb} \otimes e^{2} \otimes d^{2} \otimes \be_{\Dh}) u \\
&&+ \beta(g_{\halb}) \cdot \beta(g_{0}) \cdot \mu(e^{2} \otimes d^{2}) + \beta(g_{1}) \cdot \beta(g_{\halb}) \cdot \bar{\varrho}(e^{2} \otimes d^{2}) \\
&& + \beta(g_{1}) u^{*} \lambda(g_{\halb}^{\halb} \otimes e^{2} \otimes d^{2} \otimes \be_{\Dh}) u\\
&=& \beta(g_{0}) \cdot \mu(e^{2} \otimes d^{2}) + u^{*} \lambda(g_{\halb}^{\halb} \otimes e^{2} \otimes d^{2} 
\otimes \be_{\Dh}) u + \beta(g_{1}) \cdot \bar{\varrho}(e^{2} \otimes d^{2}) \\
& = & \gamma(e^{2} \otimes d^{2}) \, ;
\end{eqnarray*}
here, for the second equation we have used j), q),  Lemma \ref{multiplicative-domain}, (\ref{-1}), (\ref{5}), 
(\ref{6}), (\ref{7}) and (\ref{7.5}). As a consequence, all of $E \otimes \Dh$ is in the multiplicative domain 
of $\gamma$, so $\gamma$ is a $*$-homomorphism. It now follows from Theorem \ref{nonunital-intertwining} that 
$E$ in fact is $\Dh$-stable.
\end{nproof}
\en

\section{$\mathrm{K}$-theory and classification}\label{ktheory}

In this section we examine the ordered $\mathrm{K}$-theory of strongly self-absorbing
$C^*$-algebras in the UCT class and derive a number of  classification results.  In particular, we show that $\mathrm{K}_{1}$ of such algebras is always trivial and that $\mathrm{K}_{0}$ can only be $0$, $\Z$ or the $\mathrm{K}_{0}$-group of a UHF algebra of infinite type. As a consequence, we give an exhaustive list of purely infinite  strongly self-absorbing $C^{*}$-algebras. In the stably finite case, we restrict ourselves to certain inductive limit algebras; it turns out that these are either projectionless or UHF of infinite type.

\bn
\label{Kunneth}
\begin{nprop}
Let $\Dh$ be a strongly self-absorbing $C^*$-algebra satisfying the Universal
Coefficients Theorem. Then, $\mathrm{K}_1 \Dh =0$, and $\mathrm{K}_0 \Dh$
is group isomorphic to one of $0$, $\Z$, or the $\mathrm{K}_0$-group 
of a UHF algebra of infinite type. If $\mathrm{K}_{0}\Dh \cong \Z$, then $\be_{\Dh}$ represents a generator of $\mathrm{K}_{0}\Dh$.
\end{nprop}

\begin{nproof}
The UCT yields short exact sequences
\[
0 \rightarrow \mathrm{K}_* \Dh \otimes \mathrm{K}_* \Dh \rightarrow
\mathrm{K}_*(\Dh \otimes \Dh) \rightarrow \mathrm{Tor}(\mathrm{K}_* \Dh ,\mathrm{K}_* \Dh)
\rightarrow 0
\]
and
\[
0 \rightarrow \mathrm{K}_* \Dh \otimes \mathrm{K}_* \C \rightarrow
\mathrm{K}_*(\Dh \otimes \C) \rightarrow \mathrm{Tor}(\mathrm{K}_* \Dh,\mathrm{K}_* \C)
\rightarrow 0.
\]
The inclusion 
\[
\mathrm{K}_* \Dh  \otimes \mathrm{K}_* \C \rightarrow
\mathrm{K}_*(\Dh \otimes \C)
\]
is an isomorphism, since  
$\mathrm{Tor}(\mathrm{K}_* \Dh,\mathrm{K}_* \C)=0$.
Since $\Dh$ is strongly self-absorbing, the map $\mathrm{id}_{\Dh} 
\otimes \mathrm{id}_{\C} \cdot \be_{\Dh}$ induces an isomorphism
\[
\mathrm{K}_*(\mathrm{id}_{\Dh} \otimes \mathrm{id}_{\C} \cdot \be_{\Dh}):
\mathrm{K}_*(\Dh \otimes \C) \longrightarrow \mathrm{K}_*(\Dh \otimes \Dh)
\]
which, by naturality, factorises through the inclusion of the first short
exact sequence.  Said inclusion, therefore, is an isomorphism, whence
$\mathrm{Tor}(\mathrm{K}_* \Dh,\mathrm{K}_* \Dh)=0$ and
$\mathrm{K}_* \Dh$ is torsion free.  The map from $\mathrm{K}_* \Dh
\otimes \mathrm{K}_* \C$ to $\mathrm{K}_* \Dh \otimes \mathrm{K}_* \Dh$
does not meet $\mathrm{K}_1 \Dh  \otimes \mathrm{K}_1 \Dh$, and since the
composition of this map with the inclusion of $\mathrm{K}_* \Dh \otimes 
\mathrm{K}_* \Dh$ into $\mathrm{K}_*(\Dh \otimes \Dh)$ is also an isomorphism,
we have $\mathrm{K}_1 \Dh \otimes \mathrm{K}_1 \Dh =0$.  This implies that
$\mathrm{K}_1 \Dh =0$.

It remains to examine $\mathrm{K}_0 \Dh$.  
From the analysis above we have that the inclusion $\psi = \id_{\mathrm{K}_{0}} \otimes [\be_{\Dh}]$ of 
$\mathrm{K}_0 \Dh$ into $\mathrm{K}_0 \Dh \otimes
\mathrm{K}_0 \Dh$ as the first factor is an isomorphism.
Suppose that $\mathrm{K}_0 \Dh$ contains a non-zero (and necessarily 
torsion free) subgroup $H$ which is independent
of $\langle [\be_{\Dh}] \rangle$.  Then, the image
of $\psi$ will fail to meet $\mathrm{K}_0 \Dh \otimes H$ 
non-trivially, contradicting the fact that $\psi$ is an isomorphism.
Thus, every subgroup of $\mathrm{K}_0 \Dh$ meets $\langle [\be_{\Dh}] \rangle$.  
Let $x \in \mathrm{K}_0 \Dh$, and let $m,n$ be integers
such that $mx = n[\be_{\Dh}]$.  If $mx' = n[\be_{\Dh}]$ for some
$x' \in \mathrm{K}_0 \Dh$, then $m(x-x')=0$.  This implies that
$x=x'$, since $\mathrm{K}_0 \Dh$ is torsion free.  Thus, each 
element of $\mathrm{K}_0 \Dh$ is a rational multiple of $[\be_{\Dh}]$. Now 
if $[\be_{\Dh}]=0 \in \mathrm{K}_{0}\Dh$, then $\mathrm{K}_{0} \Dh$ is the 
trivial group. Otherwise,  there is an embedding $\iota$ of $\mathrm{K}_0 \Dh$ into $\Q$
which sends $[\be_{\Dh}]$ to $1$. It is straightforward to check that
\begin{equation}
\label{K_0-products}
\iota \circ \psi^{-1}(x \otimes y) = \iota(x) \cdot \iota(y)
\end{equation}
for all $x, \, y \in \mathrm{K}_{0}\Dh$. This in particular implies that 
$\iota(\mathrm{K}_{0}\Dh) \cap \Q_{+}$ cannot contain a minimal element 
other than $1$. Therefore, if $\mathrm{K}_{0}\Dh \cong \Z$, then it is 
generated by $[\be_{\Dh}]$.  The argument also shows that, if $\mathrm{K}_0 \Dh$ is not
isomorphic to $0$ or $\Z$, then it is
infinitely generated.  To complete the proof, one only has to verify that if
$1/m \in \iota(\mathrm{K}_0 \Dh)$, then so is $1/m^2$, but this follows directly 
from (\ref{K_0-products}).
\end{nproof}
\en

\bn
\label{pi-classification}
By Theorem \ref{dichotomy}, a strongly self-absorbing $C^{*}$-algebra is either purely 
infinite or stably finite with unique tracial state. In the former case, Proposition 
\ref{Kunneth} together with the Kirchberg--Phillips classification theorem (cf.\ 
\cite{R1}, Theorem 8.4.1) allows us to write down an exhaustive list -- at least of 
those algebras in the UCT class: 

\begin{ncor}
Suppose $\Dh$ is a separable purely infinite strongly self-absorbing $C^{*}$-algebra 
which satisfies the Universal Coefficients Theorem. Then $\Dh$ is either $\Oh_{2}$, 
$\Oh_{\infty}$ or a tensor product of $\Oh_{\infty}$ with a UHF algebra of infinite type. 
\end{ncor}
\en

\bn
\label{K_0-dichotomy}
In the stably finite case, the situation is less clear and we only have partial results.  
The next proposition says that the problem of classifying stably finite strongly self-absorbing 
$C^{*}$-algebras falls in two parts. 

\begin{nprop}
If $\Dh$ is a stably finite strongly self-absorbing $C^*$-algebra, then
it is either projectionless, or contains projections of arbitrarily small
trace.
\end{nprop}

\begin{nproof}
Suppose that $\Dh$ contains a non-trivial projection $p$.  Since $\Dh$ is simple, 
the unique tracial state $\tau$ on $\Dh$ is faithful, whence $0 < \tau(p) < 1$. For 
any $k \in \N$ there is an isomorphism between $\Dh$ and $\Dh^{\otimes k}$ which 
takes $\tau$ to $\tau^{\otimes k}$ (this obviously is a tracial state on $\Dh^{\otimes k}$, 
and it has to be unique). Now if $k$ is chosen large enough, $\tau^{\otimes k}(p^{\otimes k})=
\tau(p)^{k}$ becomes arbitrarily small.  
\end{nproof}
\en

\bn
There are not many classification results available for projectionless $C^{*}$-algebras. Currently, 
the most general such result applicable to our setting is the classification theorem  of \cite{M}, 
which implies that $\Zh$ is the only projectionless  strongly self-absorbing example in the class of 
simple inductive limits of circle algebras with dimension drops. We do not have any information for more general 
classes of inductive limit algebras in the projectionless case, but if nontrivial projections do 
exist we are in a much better position.   As a first step in this direction, we draw some conclusions 
about the structure of a strongly self-absorbing $C^*$-algebra when it is an 
inductive limit of recursive subhomogeneous $C^*$-algebras. 

Recall that a recursive subhomogeneous $C^*$-algebra is given by the following two part recursive
definition (Definition 1.1 of \cite{P4}):
\begin{enumerate}
\item If $X$ is a compact Hausdorff space, then $\Ch(X,M_n)$ is a recursive
subhomogeneous $C^*$-algebra for every $n \in \mathbb{N}$.
\item If $A$ is a recursive subhomogeneous $C^*$-algebra, $X$ is a compact Hausdorff space,
$X^{(0)} \subseteq X$ is closed, $\phi:A \to \Ch(X^{(0)},M_n)$ is any unital
homomorphism, and $\rho:\Ch(X,M_n) \to \Ch(X^{(0)},M_n)$
is the restriction homorphism, then the pullback
\[
A \oplus_{\Ch(X^{(0)},M_n)} \Ch(X,M_n) :=
\{(a,f) \in A \oplus \Ch(X,M_n) \, | \, \phi(a)=\rho(f)\}
\]
is a recursive subhomogeneous $C^*$-algebra.
\end{enumerate}
Given a recursive subhomogeneous $C^{*}$-algebra $A$, one can choose a so-called recursive subhomogeneous decomposition (which is highly nonunique in general). To each such decomposition one associates a compact Hausdorff space, say $Y$, called the \emph{total space}, and its \emph{topological dimension} $\dim Y$. $A$ may then be viewed as an algebra of continuous matrix-valued functions on $Y$. The \emph{topological dimension function}
of $A$ is the map $d:Y \to \mathbb{Z}^+$ which assigns to $y \in Y$ the covering dimension
of the connected component of $Y$ containing $y$. We shall only consider recursive subhomogeneous $C^{*}$-algebras which admit a decomposition with finite topological dimension. The largest matrix size of a finite-dimensional representation of $A$ is called the \emph{maximum matrix size} of $A$, while the smallest such matrix size is the \emph{minimum matrix size}. We refer the interested reader to \cite{P4}, \cite{P3} and \cite{Wi3} for more details on recursive subhomogeneous $C^*$-algebras.

Let $(A_i,\phi_{ij})$ be a direct system of recursive subhomogeneous $C^*$-algebras, where
each $A_i$ is equipped with a choice of total space $X_i$ and topological dimension function
$d_i$.  The system is said to have \emph{slow dimension growth} if for every $i$,
every projection $p \in M_{\infty}(A_i)$, and every $N \in \mathbb{N}$, there
is $j_0$ such that for all $j \geq j_0$ and $x \in X_j$ one has
\[
\phi_{ij}(p)(x) = 0 \ \ \mathrm{or} \ \ \mathrm{rank}(\phi_{ij}(p)(x)) \geq Nd_j(x).
\]
If we do not allow $\phi_{ij}(p)(x) = 0$ for $p \neq 0$, the system is said to have
\emph{strict slow dimension growth}.

We say that an inductive limit $A$ of recursive subhomogeneous $C^*$-algebras has slow
dimension growth (resp.\ strict slow dimension growth) if it can be written as the 
limit of a direct system of recursive subhomogeneous $C^*$-algebras with slow dimension
growth (resp.\ strict slow dimension growth).
\en

\bn
\label{unperf-K_0}
\begin{ntheorem}
Let $A$, $\Dh$ be unital inductive limits of recursive subhomogeneous $C^*$-algebras.  
Suppose that $\Dh$ is strongly self-absorbing, and that $A$ is $\Dh$-stable.
Then, $A$ has slow dimension growth.  If, in addition, $A$ is simple, then it
has strict slow dimension growth.  Finally, if $A$ and $\Dh$ are AH algebras,
then $A$ has very slow dimension growth in the sense of Gong (\cite{G}).
\end{ntheorem} 

\begin{nproof} 
Write $A = \lim_{i \rightarrow \infty} (A_i,\phi_i)$ and $\Dh = \lim_{i \to \infty}(\Dh_i,\gamma_i)$,
where each $A_i$ and $\Dh_i$ is a recursive subhomogeneous $C^*$-algebra with total spaces $X_i$ and $Y_i$,
respectively.  Put $a_i = \mathrm{dim}(X_i)$ and $b_i = \mathrm{dim}(Y_i)$.  Since the class of recursive
subhomogeneous $C^*$-algebras is closed under taking quotients, we may assume that the $\phi_i$ and $\gamma_i$
are injective.  Consider the commutative diagram
\[
\xymatrix{
{A \otimes \Dh}\ar[rr]^{\mathrm{id}_A \otimes \mathrm{id}_\Dh \otimes 1_\Dh}&&{A \otimes \Dh^{\otimes 2}}\ar[rrr]^{\mathrm{id}_A 
\otimes \mathrm{id}_{\Dh^{\otimes 2}} \otimes 1_\Dh}&&&{\cdots}\ar[r]&{A \otimes \Dh^{\otimes \infty}}\\
{\vdots}\ar[u]&& {\vdots}\ar[u]&&&& \\
{{A_2 \otimes \Dh_2}}\ar[rr]^{\mathrm{id}_{(A_2 \otimes \Dh_2)} \otimes \be_{\Dh_2}}\ar[u]^{\phi_2 \otimes \gamma_2}&&
{A_2 \otimes \Dh_2^{\otimes 2}}\ar[rrr]^{\mathrm{id}_{(A_2 \otimes \Dh_2^{\otimes 2})} \otimes \be_{\Dh_2}}\ar[u]^{\phi_2 \otimes \gamma_2^{\otimes 2}}&&&\cdots& \\
{A_1 \otimes \Dh_1}\ar[rr]^{\mathrm{id}_{(A_1 \otimes \Dh_1)} \otimes \be_{\Dh_1}}\ar[u]^{\phi_1 \otimes \gamma_1}&&
{A_1 \otimes \Dh_1^{\otimes 2}}\ar[rrr]^{\mathrm{id}_{(A_1 \otimes \Dh_1^{\otimes 2})} \otimes \be_{\Dh_2}}\ar[u]^{\phi_1 \otimes \gamma_1^{\otimes 2}}&&&\cdots& .
}
\]
Label the algebra $A_j \otimes \Dh^{\otimes i}_j$ with the ordered pair $(i,j)$.  Let $s((i,j),(k,l))$ 
be the path from $(i,j)$ to $(k,l)$ obtained by composing the
horizontal path from $(i,j)$ to $(k,j)$ with the vertical path from $(k,j)$ to $(k,l)$.  For any sequence
$(i_n,j_n) \in \N \times \N$ which is strictly increasing in both variables we have
\[
\lim_{n \rightarrow \infty}\left(A_{j_n} \otimes \mathcal{D}_{j_n}^{\otimes i_n},s((i_n,j_n),(i_{n+1},j_{n+1}))\right) 
\cong A \otimes \mathcal{D}^{\otimes \infty} \cong A.  
\]

The $A_{j_n} \otimes \mathcal{D}_{j_n}^{\otimes i_n}$ are recursive subhomogeneous algebras which admit recursive subhomogeneous decompositions with total spaces $X_{j_{n}} \times Y_{j_{n}}^{i_{n}}$ and topological dimension $a_{j_{n}} + i_{n} b_{j_{n}}$, cf.\ also  Proposition 3.4 of \cite{P4}.  Let 
$j_n=n$, and let $(\varepsilon_n)$ be a sequence of positive tolerances converging to zero. 
Set 
\[
s_n := s((i_n,n),(i_{n+1},n+1)), 
\]
and assume that for $k < n$, $i_k$ has been chosen
with the following property:  if $p \in A_{k-1} \otimes \mathcal{D}_{k-1}^{\otimes i_{k-1}}$ is a projection,
then
\[
\mathrm{rank}(s_{k-1}(p)(x)) \geq \frac{1}{\varepsilon_k}(a_k + i_k b_k)
\]
for every $x$ in the total space of $A_k \otimes \Dh_k^{\otimes i_k}$ such that $s_{k-1}(p)(x) \neq 0$.
We prove that $i_n$ can be chosen in a like manner.    

Lemma 1.8 of \cite{P3} states that if $B = \lim_{i \to \infty}(B_i,\eta_i)$ 
is a simple, unital, and infinite-dimensional inductive limit of
recursive subhomogeneous $C^*$-algebras with injective and unital connecting morphisms (in particular,
we could take $B = \Dh$), then for any projection $p \in B_i$ and $N \in \mathbb{N}$ 
there exists $j_0$ such that for every $j \geq j_0$ one has
\[
\mathrm{rank}(\eta_{ij}(p))(x) \geq N
\]
for every $x$ in the total space of $B_j$.  In particular, the rank of the $\be_{\Dh_i}$ may
be assumed to be greater than two for every $i \in \mathbb{N}$.  

Choose $i_n$ to satisfy
\[
\frac{2^{(i_n - i_{n-1})}}{a_n + i_n b_n} > \frac{1}{\varepsilon_n},
\] 
and let $p \in A_{n-1} \otimes \Dh_{n-1}^{\otimes i_{n-1}}$ be a projection.  Then,
$s_{n-1}(p)$ may be viewed as an elementary tensor $q \otimes \be_{\Dh_n}^{\otimes i_n - i_{n-1}}$,
where $q \in A_n \otimes \Dh_n^{\otimes i_{n-1}}$ is a projection.  Since 
\[
\mathrm{rank}(\be_{\Dh_n}^{\otimes i_n - i_{n-1}}) \geq \frac{1}{\varepsilon_n}(a_n + i_n b_n)
\]
at every point in the total space of $\Dh_n^{\otimes i_n - i_{n-1}}$ (which we can choose to be homeomorphic
to $(Y_n)^{i_n-i_{n-1}}$ --- cf.\ Proposition 3.4 of \cite{P4}), the same is true of
\[
\mathrm{rank}(q \otimes \be_{\Dh_n}^{\otimes i_n - i_{n-1}})
\]
over each point in the total space of $A_n \otimes \Dh_n^{\otimes i_n}$, as required.  Thus,
$A$ has slow dimension growth.

If $A$ is simple, then one may apply Lemma 1.8 of \cite{P3} to conclude that the projection $p$
above may be chosen to have non-zero rank over every point in the total space of $A_{n-1} \otimes 
\Dh_{n-1}^{\otimes i_{n-1}}$.  The rank estimates above then apply over every point in the total
space of $A_n \otimes \Dh_n^{\otimes i_n}$, and $A$ has strict slow dimension growth.

Finally, suppose that $A$ and $\Dh$ are AH algebras.  It follows from \cite{EGL} that we may
assume the connecting morphisms in the inductive sequences for $A$ and $\Dh$ to be unital and
injective, whence the closure of the class of homogeneous $C^*$-algebras under tensor products
allows us to repeat the above proof inside the class of AH algebras.  It follows that $A$ has
very slow dimension growth in the sense of \cite{G}.
\end{nproof}
\en

\bn
The principal theorems of \cite{P3} now yield:

\begin{ncor}
Let $A$ and $\Dh$ be as in Theorem \ref{unperf-K_0}.  Then,
\begin{enumerate}
\item[(i)] The map $\mathcal{U}(A)/\mathcal{U}(A)_0 \to \mathrm{K}_1 A$ is an isomorphism;
\item[(ii)] if $A$ is simple, then $\mathrm{K}_0 A$ is weakly unperforated;
\item[(iii)] if $A$ is simple, then the projections in $M_{\infty}(\mathcal{D})$ satisfy cancellation;
\item[(iv)] if $A$ is simple, then it satisfies Blackadar's second fundamental
comparability property. 
\end{enumerate}
\end{ncor}
\en

\bn
In \cite{EGL} it is shown that  simple
unital AH algebras with very slow dimension growth are classified by their Elliott invariants; it is also known that such algebras contain nontrivial projections.
This leads us to the following corollary of \ref{Kunneth}, \ref{K_0-dichotomy}   and \ref{unperf-K_0}:

\label{AH-classified}
\begin{ncor}
The strongly self-absorbing AH $C^*$-algebras are classified by the Elliott invariant; they are precisely the UHF algebras of infinite type.
\end{ncor}
\en

\bn
\label{weakdiv-appdiv}
Recall that a simple unital $C^{*}$-algebra $A$ is approximately divisible if it admits an approximately 
central sequence of unital $*$-homomorphisms of $M_{2} \oplus M_{3}$ into $A$. It is weakly divisible, if 
for each projection $p$ in $A$ there is a unital $*$-homomorphism of $M_{2} \oplus M_{3}$ into $pAp$. An 
approximately divisible $C^{*}$-algebra has real rank zero if and only if projections in $A$ separate traces; 
it satisfies Blackadar's second fundamental comparability property (cf.\ \cite{R1} and the references therein) 
and it is $\Zh$-stable by \cite{TW2}. If $A$ has real rank zero, it is weakly divisible by \cite{RP}. The 
next observation says that all these properties coincide for strongly self-absorbing $C^{*}$-algebras with projections and, in the purely infinite case, are automatically fulfilled. 

\begin{nprop}
Consider the following conditions for a strongly self-absorbing $C^{*}$-algebra $\Dh$:
\begin{itemize}
\item[(i)] $\Dh$ contains a nontrivial projection and satisfies Blackadar's second fundamental comparability property.
\item[(ii)] $\Dh$ is weakly divisible.
\item[(iii)] $\Dh$ is approximately divisible.
\item[(iv)] $\Dh$ has real rank zero.
\item[(v)]  $\Dh$ is $\Zh$-stable and contains a nontrivial projection.
\end{itemize}
These conditions are all equivalent; they are satisfied if $\Dh$ is purely infinite.  
\end{nprop}

\begin{nproof}
First, assume $\Dh$ to be stably finite.\\
(i) $\Rightarrow$ (ii): Let $p \in \Dh$ be a projection. Identifying $\Dh$ with $\Dh^{\otimes \infty}$, $p$ is close, hence Murray--von Neumann equivalent, to a projection of the form $q \otimes \be_{\Dh^{\otimes \infty}}$ for some $k \in \N$ and a projection $q \in \Dh^{\otimes k}$. But then $p \Dh p $ is isomorphic to $(q \otimes \be_{\Dh^{\otimes \infty}}) (\Dh^{\otimes k} \otimes \Dh^{\otimes \infty}) (q \otimes \be_{\Dh^{\otimes \infty}}) = (q \Dh^{\otimes k} q) \otimes \Dh^{\otimes \infty}$, whence $p \Dh p$ is $\Dh$-stable. Therefore, it suffices to map $M_{2} \oplus M_{3}$ unitally to $\Dh$. Since $\Dh$ contains projections arbitrarily small in trace,  
using comparison one finds a projection $q \in \Dh$ such that $1/3<\tau(q)<1/2$ (by adding up sufficiently 
many projections which are all equivalent, pairwise orthogonal and small in trace). Again by comparison, 
$q$ is equivalent to a subprojection of $\be_{\Dh}-q$; this defines a (nonunital) embedding $\alpha: M_{2} 
\to \Dh$ with $\alpha(e_{11})=q$. Set $q':= \be_{\Dh} - \alpha(\be_{M_{2}})$, then $\tau(q')< 1/3<\tau(q)$ 
and (using comparison oncemore) there is another embedding $\beta:M_{2} \to \Dh$, this time with $\beta(e_{11})=q'$ 
and $\beta(e_{22}) \le q$. But now it is straightforward to check that the $C^{*}$-subalgebra of $\Dh$ 
generated by $\alpha(M_{2})$ and $\beta(M_{2})$ is in fact isomorphic to $M_{2} \oplus M_{3}$.\\
(ii) $\Rightarrow$ (iii): Map $M_{2} \oplus M_{3}$ unitally to $\Dh$ and take an approximately central sequence 
of unital embeddings of $\Dh$ into itself; the compositions will yield an approximately central sequence of 
unital $*$-homomorphisms from $M_{2} \oplus M_{3}$ to $\Dh$.\\
(iii) $\Rightarrow$ (iv): An approximately divisible $C^{*}$-algebra has real rank zero if and only if 
projections separate tracial states -- which is always true in the unique trace case.\\
(iv) $\Rightarrow$ (ii) was shown in \cite{RP}.\\
(iii) $\Rightarrow$ (v): A unital approximately divisible $C^{*}$-algebra clearly contains a nontrivial projection. 
$\Zh$-stability essentially follows from Theorem \ref{rordams-intertwining} and the special inductive limit structure 
of $\Zh$; a complete proof can be found in \cite{TW2}.\\
(v) $\Rightarrow$ (i) follows from \cite{J}.

If $\Dh$ is purely infinite, it has real rank zero, so (iv) holds. Condition (i) is trivially true, since $\Dh$ admits no tracial states. The proofs of the implications (iv)  $\Rightarrow$ (ii) $\Rightarrow$ (iii) $\Rightarrow$ (v) apply verbatim. 

We wish to point out that implications (iv) $\Rightarrow$ (iii) and (iv) $\Rightarrow$ (v) can also be deduced from \cite{ElR}, Corollary 2.4.
\end{nproof}
\en
 
\bn
\label{typeI-UHF}
The above proposition has a very satisfying corollary via
recent work of N.\ Brown (\cite{B}).

\begin{ncor}
Let $\mathcal{D}$ be a strongly self-absorbing $C^*$-algebra satisfying one of the conditions of the preceding proposition. 
Suppose further that $\mathcal{D}$ is an inductive limit of type I $C^*$-algebras.
Then, $\mathcal{D}$ is a UHF algebra of infinite type.
\end{ncor}

\begin{nproof}
Since $\mathcal{D}$ is strongly self-absorbing, it is unital, simple, and nuclear.  The hypothesis
that $\mathcal{D}$ is an inductive limit of type I $C^*$-algebras implies that $\mathcal{D}$ satisfies
the UCT, and that its unique trace satisfies Definition 6.1 of \cite{B}.  The approximate
divisibility of $\mathcal{D}$ implies that $\mathcal{D}$ has stable rank one
and weakly unperforated $\mathrm{K}$-theory.  Since $\Dh$ has real rank zero and has weakly unperforated $\mathrm{K}_{0}$-group, we may apply Corollary 7.9 of \cite{B}, which 
implies that $\mathcal{D}$ is tracially AF.  Since $\mathcal{D}$ contains nontrivial projections, its $\mathrm{K}_0$-group
must be that of a UHF algebra of infinite type by Propositions \ref{Kunneth} and \ref{K_0-dichotomy}.  Lin's classification theorem for tracially AF algebras
(\cite{Li1}) now implies that $\mathcal{D}$ is UHF of infinite type.
\end{nproof}
\en

\bn
\label{RSH-UHF}
We note a variation of the above to point out a slightly different point of view. 

\begin{ncor}
Let $\Dh$ be a strongly self-absorbing inductive limit of recursive subhomogeneous algebras. Then, $\Dh$ is either projectionless or a UHF algebra of infinite type.
\end{ncor}

\begin{nproof}
Since $\Dh$ has strict slow dimension growth by Theorem \ref{unperf-K_0}, it satisfies 
Blackadar's second fundamental comparability property by \cite{P3}. Now if $\Dh$ contains 
a nontrivial projection, then condition (i) of Proposition \ref{weakdiv-appdiv} holds and 
$\Dh$ is UHF of infinite type by Corollary \ref{typeI-UHF}. 
\end{nproof}
\en 

\bn
In view of the preceding discussion, it seems natural to ask the following 

\begin{nquestion}
Are there any other stably finite strongly self-absorbing $C^{*}$-algebras than $\Zh$ and the UHF algebras of infinite type? Are these, at least, the only examples which are limits of (recursive) subhomogeneous algebras?
\end{nquestion}
\en

\bn
\label{mutual-embeddings}
To solve the above question, the following observation might be useful. It says that the problem of classifying strongly self-absorbing $C^*$-algebras
is reduced to the problem of determining when one such algebra can be 
embedded unitally in another.

\begin{nprop}
Let $\Dh$, $\Eh$, be strongly self-absorbing $C^*$-algebras, and suppose that
there exist unital embeddings $\iota_\Dh:\Dh \rightarrow \Eh$ and $\iota_\Eh:\Eh \rightarrow \Dh$.
Then, $\Dh$ and $\Eh$ are isomorphic.
\end{nprop}

\begin{nproof}
By Proposition \ref{t-self-absorbing}, $\iota_\Dh$ gives rise to
an approximately central sequence of unital embeddings of $\Dh$ into
$\Eh$, whence $\Eh$ is $\Dh$-stable by Theorem \ref{rordams-intertwining}.  Similarly,
$\Dh \otimes \Eh \cong \Dh$.
\end{nproof}
\en

\bn
We close with some remarks on topological covering dimension of $\Dh$-stable $C^{*}$-algebras. 
Recall from \cite{KW} and \cite{Wi4}, that the known examples of strongly self-absorbing $C^{*}$-algebras 
have topological dimension (i.e., decomposition rank or dimension as an AH or ASH algebra, respectively) 0, 1 or infinity. \\
If $A$ is a simple AH algebra with slow dimension growth in the AH sense, then $A \otimes \Zh$, which is a $\Zh$-stable 
ASH algebra {\it a priori}, obviously has  slow dimension growth in the ASH sense. However, in the case where $A \otimes \Zh$ has real rank zero, it will follow from results in \cite{Wi6} and \cite{Li1} that $A \otimes \Zh$ is isomorphic to $A$ (hence is AH) and, in fact, has bounded topological dimension. It is tempting to ask whether a similar reasoning (which might  involve some type of reduction step as, for example, in \cite{D1}) works in a more general setting. This is interesting for strongly self-absorbing $C^{*}$-algebras as well as for $\Dh$-stable $C^{*}$-algebras. In view of the classification results of \cite{Wi4}, \cite{Wi5} and \cite{Wi6} it seems also natural to rephrase the question in terms of the decomposition rank (cf.\ \cite{KW} and \cite{Wi3}):

\begin{nquestion}
Does every separable, strongly self-absorbing limit $\Dh$ of recursive subhomogeneous algebras have
bounded topological dimension or, at least, finite decomposition rank? Does the respective statement hold for 
$\Dh$-stable simple limits of recursive subhomogeneous algebras?
\end{nquestion}
\en




\end{document}